\documentclass[reqno]{amsart}
\usepackage{hyperref}
\usepackage{setspace}
\usepackage{graphicx}
\usepackage{amssymb}
\usepackage{mathrsfs}

\begin{document}
    \title[\hfilneg ]
    { Existence and continuation of solutions of Hilfer-Katugampola-type fractional
differential equations }

 \author[\hfil\hfilneg]{Ahmad Y. A. Salamooni, D. D. Pawar }
 \address{Ahmad Y. A. Salamooni \newline
     School of Mathematical Sciences, Swami Ramanand Teerth Marathwada University, Nanded-431606, India}
      \email{ayousss83@gmail.com}

 \address{D. D. Pawar \newline
    School of Mathematical Sciences, Swami Ramanand Teerth Marathwada University, Nanded-431606, India}
     \email{dypawar@yahoo.com}


    \begin{abstract}
    This article contains a new discussion for the generalized fractional
Cauchy-type problem involving Hilfer-Katugampola-type fractional derivative.
We study an existence and continuation of its solution. Firstly, we establish a new theorems of local existence.
Then, we deduce a continuation theorems for a general fractional differential equations.
By applying continuation theorems deduced in this article, we present several global existence results. 
Moreover, the examples are given to illustrate our main results.
      \\\\ \textbf{AMS Classification- 34A08; 26A33; 34A12.}
    \end{abstract}

    \maketitle \numberwithin{equation}{section}
    \newtheorem{theorem}{Theorem}[section]
    \newtheorem{lemma}[theorem]{Lemma}
    \newtheorem{definition}[theorem]{Definition}
    \newtheorem{example}[theorem]{Example}

    \newtheorem{remark}[theorem]{Remark}
    \allowdisplaybreaks
  \section{Introduction}
During last years, the operators  number of fractional integration and differentiation has
been increasing with new definitions for them, which correspond to  the Riemann-Liouville, Caputo,  Hadamard, Hilfer, etc.[1-2].
Actually, there exists more than one new definition for fractional derivatives,
which are Katugampola, Caputo-Katugampola and  Hilfer-Katugampola for more details see [3-5].
\par Here, we mention that the authors and researchers have interest to solve fractional differential equations,
where they have studied the existence and uniqueness theorems of a solution for fractional differential equations
by applying a fixed point theory on a finite interval $[0,T]$ and obtained a global existence
of solutions by establishing a local existence theorems and a continuation theorems on the half axis $[0,+\infty),$ [6-19] and references therein. \par C. Kau
et al. [20], they have obtained the existence and continuation theorems for the following Riemann-Liouville type fractional differential equations
\begin{equation}
\left\{\begin{matrix} _{RL}D^{\alpha}_{0,t}x(t)=\varphi(t,x(t)),\quad\quad t\in(0,+\infty),~~0<\alpha<1,
\\\\_{RL}I^{1-\alpha}_{0+}x(t)\big|_{t=0}=x_{0},\quad\quad\quad
\quad\quad\quad\quad\quad\quad\quad\end{matrix}\right.
\end{equation}
where $_{RL}D^{\alpha}_{0,t}$ is the Riemann-Liouville-type
fractional derivative of order $\alpha.$
\\Moreover, S. P. Bhairat [21], found the existence and continuation of solutions for the following Hilfer fractional
differential equations
\begin{equation}
\left\{\begin{matrix} D^{\alpha,\beta}_{0+}x(t)=\varphi(t,x(t)),\quad\quad t\in(0,+\infty),~~0<\alpha<1,0\leq\beta\leq1,
\\\\I^{1-\gamma}_{0+}x(t)\big|_{t=0}=x_{0},\quad\gamma=\alpha+\beta(1-\alpha),\quad\quad\quad
\quad\quad\quad\quad\quad\quad\quad\end{matrix}\right.
\end{equation}
where $D^{\alpha,\beta}_{0+}$ is the Hilfer-type
fractional derivative of order $\alpha$ and type $\beta.$
\par The aim of this article is to develop the existence and uniqueness theory. Firstly,  we establish the local existences of Hilfer-Katugampola and the system of Hilfer-Katugampola fractional differential equations, then we study the continuation theorems of Hilfer-Katugampola fractional differential equations to extend the existence of global solutions.
\par In this paper, we consider the Cauchy-type problem involving Hilfer-Katugampola-type
fractional derivative with initial value problems
\begin{equation}
\left\{\begin{matrix} _{\rho}D^{\alpha,\beta}_{0+}x(t)=\varphi(t,x(t)),\quad\quad t\in(0,+\infty),~~0<\alpha<1,0\leq\beta\leq1,
\\\\_{\rho}I^{1-\gamma}_{0+}x(t)\big|_{t=0}=x_{0},\quad\gamma=\alpha+\beta(1-\alpha),\quad\quad\quad
\quad\quad\quad\quad\quad\quad\quad\end{matrix}\right.
\end{equation}
where $\rho>0~and~_{\rho}D^{\alpha,\beta}_{0+}$ is the Hilfer-Katugampola-type
fractional derivative of order $\alpha$ and type $\beta,$ [5] and
$\varphi:\mathbb{R^{+}}\times\mathbb{R}\rightarrow\mathbb{R}$ has a weak singularity with respect to $t$
and here $\varphi$ satisfies a Lipschitz condition
\[|\varphi(t,x(t))-\varphi(t, y(t))|\leq A|x(t)- y(t)|\] where $A>0$ is Lipschitz constant.
\par And we consider a system of fractional differential equations with general initial value problems
\begin{equation}
\left\{\begin{matrix}_{\rho}D^{\alpha,\beta}_{0+}x_{1}(t)=\varphi_{1}(t,x_{1}(t),x_{2}(t),...,x_{n}(t)),
\quad\quad\quad\quad\quad\quad\quad\quad\quad\\\\
_{\rho}D^{\alpha,\beta}_{0+}x_{2}(t)=\varphi_{2}(t,x_{1}(t),x_{2}(t),...,x_{n}(t)),
\quad\quad\quad\quad\quad\quad\quad\quad\quad\\ \vdots\\
_{\rho}D^{\alpha,\beta}_{0+}x_{n}(t)=\varphi_{n}(t,x_{1}(t),x_{2}(t),...,x_{n}(t)),
\quad\quad\quad\quad\quad\quad\quad\quad\quad\\\\
_{\rho}I^{1-\gamma}_{0+}x_{\ell}(t)\big|_{t=0}=x_{0}\quad
\gamma=\alpha+\beta(1-\alpha),\quad\ell=1,2,...,n,
\quad\quad\quad\end{matrix}\right.
\end{equation}
where $0<\alpha<1,0\leq\beta\leq1,~\gamma=\alpha+\beta(1-\alpha),\rho>0~and~\varphi_{\ell}:\mathbb{R^{+}}\times\mathbb{R}^{n}\rightarrow\mathbb{R}$
also has a weak singularity with respect to $t$ and here $\varphi_{n}(t,x_{1}(t),x_{2}(t),...,x_{n}(t))$ satisfy a Lipschitz condition
\[|\varphi_{j}(t,x_{1}(t),x_{2}(t),...,x_{n}(t))-
\varphi_{j}(t, y_{1}(t),y_{2}(t),...,y_{n}(t))|\leq \sum_{j=1}^{n}A_{j}|x_{j}(t)- y_{j}(t)|\]
 where $A_{j}>0,~j=1,2,...,n$ is Lipschitz constants.
\par The remaining parts of this paper is ordered as below: \\In section 2, we present some basic notations,
definitions and lemmas used in our main results. Section 3, includes the study of a local existence of solutions,
in which we obtain the new local existence theorems for the initial value problems $(1.3) ~~and~~ (1.4).$
Two continuation theorems with global existence theorems for the initial value problems $(1.3)$ are given in Section 4. The last section contains concluding remarks. 

\section{Preliminaries}
In this section, we introduce some notations, definitions and lemmas from theory of fractional calculus which will be used later.
\\\textbf{\ Definition 2.1.[3]} Let $\Omega=[0,T]$ is a finite interval and $\rho>0,$ the weighted space $C_{1-\gamma,\rho}[0,T]$ of
continuous functions $\varphi$ on $(0,T]$ is defined by
\[C_{1-\gamma,\rho}[0,T]=\{\varphi:(0,T]\rightarrow\mathbb{R}:[(t^{\rho}/\rho)]^{1-\gamma}\varphi(t)\in C[0,T]\}\]
with the norm
\[\|\varphi\|_{C_{1-\gamma,\rho}}=\bigg\|[ (t^{\rho}/\rho)]^{1-\gamma}\varphi(t)\bigg\|_{C},\quad C_{0,\rho}[0,T]=C[0,T].\]
The space $C_{1-\gamma,\rho}[0,T]$ is the complete metric space defined with the distance $d$ as
\[d(x_{1},x_{2})=\|x_{1}-x_{2}\|_{C_{1-\gamma,\rho}}[0,T]:=\max_{t\in[0,T]}\bigg|[(t^{\rho}/\rho)]^{1-\gamma}\big[x_{1}(t)-x_{2}(t)\big]\bigg|.\]
\\\textbf{\ Definition 2.2.[3]} Let $\Omega=(0,T]$ and $\phi:(0,\infty)\rightarrow\mathbb{R},$ the Katugampola fractional integrals $_{\rho}I_{0+}^{\alpha}\varphi$ of order $\gamma\in\mathbb{C}(\mathfrak{R}(\alpha)>0)$
is defined for $\rho>0$ as
\begin{align}
(_{\rho}I_{0+}^{\alpha}\varphi)(t)=\frac{\rho^{1-\alpha}}{\Gamma(\alpha)}
\int_{0}^{t}\frac{\tau^{\rho-1}\varphi(\tau)}{(t^{\rho}-\tau^{\rho})^{1-\alpha}}d\tau,\quad(t>0),\label{e2.1}
\end{align}
and the corresponding Katugampola fractional derivative $_{\rho}D_{0+}^{\alpha}\varphi$
is defined as
\begin{align}
(_{\rho}D_{0+}^{\alpha}\varphi)(t)&:=\big(t^{1-\rho}\frac{d}{dt}\big)^{n}
\big(_{\rho}I_{0+}^{n-\alpha}\varphi\big)(t)\nonumber\\&
=\frac{\rho^{\alpha-n+1}}{\Gamma(n-\alpha)}\big(t^{1-\rho}\frac{d}{dt}\big)^{n}
\int_{0}^{t}\frac{\tau^{\rho-1}\varphi(\tau)}{(t^{\rho}-\tau^{\rho})^{\alpha-n+1}}d\tau,\quad(t>0),\label{e2.2}
\end{align}
\\\textbf{\ Definition 2.3. [5]}  Let$0<\alpha<1,~~0\leq\beta\leq 1,~~\varphi\in C_{1-\gamma,\rho}[0,T].$ The
Hilfer-Katugampola fractional derivative $_{\rho}D^{\alpha,\beta}$ of
order $\alpha$ and type $\beta$ of $\varphi$ is defined
as
\begin{align}
(~_{\rho}D^{\alpha,\beta}\varphi)(t)&=\bigg(~_{\rho}I^{\beta(1-\alpha)}\big(t^{1-\rho}\frac{d}{dt}\big)~_{\rho}I^{(1-\alpha)(1-\beta)}\varphi\bigg)(t)
\nonumber\\&=\bigg(~_{\rho}I^{\beta(1-\alpha)}(\delta_{\rho})~_{\rho}I^{(1-\alpha)(1-\beta)}\varphi\bigg)(t)
;\quad\gamma=\alpha+\beta(1-\alpha).
\end{align}
Where $_{\rho}I^{(.)}$ is the Katugampola fractional
integral defined in (2.1).
\\\textbf{\ Lemma 2.1. [6]} Let $a<b<c,~0\leq\nu<1,~x\in C_{\nu}[a,b],~y\in C[b,c]$ and $x(b)=y(b).$ Define
\begin{equation}
z(t)=\left\{\begin{matrix} x(t)\quad\quad if \quad t\in(a,b],
\\\\y(t)\quad\quad if \quad t\in[b,c].\end{matrix}\right.
\end{equation}
Then, $z\in C_{\nu}[a,c].$
\\\textbf{\ Lemma 2.2. [22]}(Schauder fixed point Theorem) Let $U$ be a closed bounded convex subset of a Banach space $E$
and Suppose that $T:U\rightarrow U$ is completely continuous operator.
Then, $T$ has a fixed point in $U$.
\\\textbf{\ Lemma 2.3. [5]} Let $\Omega=[0,T]$ is a finite interval, $\alpha>0$ and $0\leq\nu<1.$
\\(a) If $\nu>\alpha,$ then the fractional integration operator $_{\rho}I_{0+}^{\alpha}$ is bounded
from $C_{\nu,\rho}[0,T]$ into $C_{\nu-\alpha,\rho}[0,T].$
\\(b) If $\nu\leq\alpha,$ then the fractional integration operator $_{\rho}I_{0+}^{\alpha}$ is bounded
from $C_{\nu,\rho}[0,T]$ into $C[0,T].$
\\\textbf{\ Lemma 2.4. [5]}  Let $0<\alpha<1,0\leq\beta\leq1,\gamma=\alpha+\beta(1-\alpha),$
and assume that $\varphi(t,x(t))\in C_{1-\gamma,\rho}[0,T]$ where
$\varphi:(0,T]\times\mathbb{R}\rightarrow\mathbb{R}$ be a function
for any $x\in C_{1-\gamma,\rho}[0,T].$
If $x\in C^{\gamma}_{1-\gamma,\rho}[0,T],$ then $x$ satisfies $(1.3)$
if, and only if, $x$ satisfies the second kind Volterra fractional integral equation
\begin{equation}
x(t)=\frac{x_{0}}{\Gamma(\gamma)}\big(t^{\rho}/\rho\big)^{\gamma-1}+\frac{\rho^{1-\alpha}}{\Gamma(\alpha)}
\int_{0}^{t}\frac{\tau^{\rho-1}\varphi(\tau,x(\tau))}{(t^{\rho}-\tau^{\rho})^{1-\alpha}}d\tau,\quad(t>0).
\end{equation}
\par In the light of the Lemma 2.3 (see[20]), we have the following Lemma
\\\textbf{\ Lemma 2.5.} Let $\mathcal{A}$ be the subset of $C_{1-\gamma,\rho}[0,T]$.
Then, $\mathcal{A}$ is precompact if, and only if, the following
conditions are satisfied:
\\(1) $\{\big(t^{\rho}/\rho\big)^{1-\gamma}x(t):x\in\mathcal{A}\}$ is uniformly bounded,
\\(2) $\{\big(t^{\rho}/\rho\big)^{1-\gamma}x(t):x\in\mathcal{A}\}$ is equicontinuous on $[0,T]$.

\section{The Local Existence}
In this section, we study the local existence of solutions for the initial value problems $(1.3) ~~and~~ (1.4).$
Assume that $\varphi(t,x(t))$ in $(1.3)$ and $\varphi_{\ell}(t,x_{\ell}(t)),(\ell=1,2...,n)$ in $(1.4)$ have some weak singularity with
respect to $t$ respectively. By using Schauder fixed point theorem, we have obtained new local
existence theorems.
\par For convenience, we create the following two hypothesis.
\\$(\mathcal{H}_{1})$ Assume that $\varphi:\mathbb{R^{+}}\times\mathbb{R}\rightarrow\mathbb{R}$ in $(1.3)$ is the continuous function and
there exists a constant $0\leq\lambda<1$ such that $(\mathcal{M}x)(t) = t^{\lambda}\varphi(t,x(t))$ be the continuous bounded
map from $C_{1-\gamma,\rho}[0,T]$ into $C[0,T],$ where $T$ be a positive constant.
\\$(\mathcal{H}_{2})$ Assume that $\varphi_{\ell}:\mathbb{R^{+}}\times\mathbb{R}^{n}\rightarrow\mathbb{R}$ in $(1.4)$ is the continuous function and
there exists a constant $0\leq\lambda_{\ell}<1$ such that
$(\mathcal{M}_{\ell}x_{\ell})(t) = t^{\lambda_{\ell}}\varphi_{\ell}(t,x_{1}(t),x_{2}(t),...,x_{n}(t)),(\ell=1,2...,n)$ are continuous bounded
maps from $C_{1-\gamma,\rho}[0,T]$ into $C[0,T],$ where $T$ be a positive constant.
\\\textbf{\ Theorem 3.1.} Assume that a condition $(\mathcal{H}_{1})$  is satisfied. Then the initial value problems $(1.3)$ has at least one solution
$x\in C_{1-\gamma,\rho}[0,h]$ for some $(T\geq)h>0.$
\\ \textbf{\ Proof.} Let
\begin{equation}
D=\bigg\{x\in C_{1-\gamma,\rho}[0,T]:\bigg\|x-\frac{x_{0}}{\Gamma(\gamma)}\big(t^{\rho}/\rho\big)^{\gamma-1}\bigg\|_{C_{1-\gamma,\rho}[0,T]}
=\sup_{0\leq t\leq T}\bigg|\big(t^{\rho}/\rho\big)^{1-\gamma}x(t)-\frac{x_{0}}{\Gamma(\gamma)}\bigg|\leq k\bigg\},\label{e3.1}
\end{equation}
where $k>0$ is a constant. Since the operator $\mathcal{M}$ is bounded, there exists a constant $L>0$ such that
\[\sup\big\{\big|(\mathcal{M}x)(t)\big|:t\in[0,T],~x\in D\big\}\leq L.\]
Again, let
\begin{equation}
E_{h}=\bigg\{x:x\in C_{1-\gamma,\rho}[0,T],
\sup_{0\leq t\leq T}\bigg|\big(t^{\rho}/\rho\big)^{1-\gamma}x(t)-\frac{x_{0}}{\Gamma(\gamma)}\bigg|\leq k\bigg\},\label{e3.2}
\end{equation}
where $h=\min\bigg\{\big(\frac{k\rho^{\alpha-\gamma+1}\Gamma(\alpha-\lambda+1)}{L\Gamma(1-\lambda)}\big)^{\frac{1}{\rho(\alpha-\gamma-\lambda+1)}},T\bigg\}.$ Obviously, $E_{h}\subseteq C_{1-\gamma,\rho}[0,T]$ be a nonempty, bounded closed and convex subset.
\par Note that $h\leq T$, we can regard $E_{h}$ and $C_{1-\gamma,\rho}[0,T]$ as the restrictions of $D$ and $C_{1-\gamma,\rho}[0,T]$, respectively.
Define the operator $\mathcal{N}$ as follows
\begin{equation}
(\mathcal{N}x)(t)=\frac{x_{0}}{\Gamma(\gamma)}\big(t^{\rho}/\rho\big)^{\gamma-1}+\frac{\rho^{1-\alpha}}{\Gamma(\alpha)}
\int_{0}^{t}\tau^{\rho-1}(t^{\rho}-\tau^{\rho})^{\alpha-1}\varphi(\tau)d\tau,\quad t\in[0,h].\label{e3.3}
\end{equation}
Observe that from $(\mathcal{H}_{1})$ and Lemma 2.3 we have $\mathcal{N}(C_{1-\gamma,\rho}[0,h])\subset C_{1-\gamma,\rho}[0,h].$
\par By relation $(3.3),$ for any $x\in C_{1-\gamma,\rho}[0,h],$ we obtain
\begin{align*}
\bigg|\big(t^{\rho}/\rho\big)^{1-\gamma}(\mathcal{N}x)(t)-\frac{x_{0}}{\Gamma(\gamma)}\bigg|&=
\bigg|\big(t^{\rho}/\rho\big)^{1-\gamma}\frac{\rho^{1-\alpha}}{\Gamma(\alpha)}
\int_{0}^{t}\tau^{\rho-\lambda-1}(t^{\rho}-\tau^{\rho})^{\alpha-1}\big[\tau^{\lambda}\varphi(\tau,x(\tau))\big]d\tau\bigg|\\& \leq
\frac{\big(t^{\rho}/\rho\big)^{1-\gamma}\rho^{1-\alpha}L}{\Gamma(\alpha)}
\int_{0}^{t}\tau^{\rho-\lambda-1}(t^{\rho}-\tau^{\rho})^{\alpha-1}d\tau\\&\leq
\frac{L h^{\rho(\alpha-\gamma-\lambda+1)}\Gamma(1-\lambda)}{\rho^{\alpha-\gamma+1}\Gamma(\alpha-\lambda+1)}\leq k,
\end{align*}
which yields that $\mathcal{N}E_{h}\subset E_{h}.$
\par Next, we will show that  $\mathcal{N}$ is continuous. For that let $x_{n},x\in E_{h}, \|x_{n}-x\|_{C_{1-\gamma,\rho}[0,h]}\rightarrow0$ as $n\rightarrow+\infty.$ In the light of a continuity of $\mathcal{M},$ we have $\|\mathcal{M}x_{n}-\mathcal{M}x\|_{C_{1-\gamma,\rho}[0,h]}\rightarrow0$ as $n\rightarrow+\infty.$
\\Now, Noticing that
\begin{align*}
\bigg|\big(t^{\rho}/\rho\big)^{1-\gamma}(\mathcal{N}x_{n})(t)-&\big(t^{\rho}/\rho\big)^{1-\gamma}(\mathcal{N}x)(t)\bigg|\\&=
\bigg|\big(t^{\rho}/\rho\big)^{1-\gamma}\frac{\rho^{1-\alpha}}{\Gamma(\alpha)}
\int_{0}^{t}\tau^{\rho-1}(t^{\rho}-\tau^{\rho})^{\alpha-1}\varphi(\tau,x_{n}(\tau))d\tau\\&\quad\quad\quad\quad-
\big(t^{\rho}/\rho\big)^{1-\gamma}\frac{\rho^{1-\alpha}}{\Gamma(\alpha)}
\int_{0}^{t}\tau^{\rho-1}(t^{\rho}-\tau^{\rho})^{\alpha-1}\varphi(\tau,x(\tau))d\tau\bigg|\\& \leq
\frac{\big(t^{\rho}/\rho\big)^{1-\gamma}\rho^{1-\alpha}}{\Gamma(\alpha)}
\int_{0}^{t}\tau^{\rho-\lambda-1}(t^{\rho}-\tau^{\rho})^{\alpha-1}
\bigg|\tau^{\lambda}\big[\varphi(\tau,x_{n}(\tau))-\varphi(\tau,x(\tau))\big]\bigg|d\tau\\&\leq
\frac{\big(t^{\rho}/\rho\big)^{1-\gamma}\rho^{1-\alpha}}{\Gamma(\alpha)}
\int_{0}^{t}\tau^{\rho-\lambda-1}(t^{\rho}-\tau^{\rho})^{\alpha-1}d\tau
\big\|\mathcal{M}x_{n}-\mathcal{M}x\big\|_{[0,h]}.
\end{align*}
 Then, we have
\begin{equation}
\big\|\mathcal{N}x_{n}-\mathcal{N}x\big\|_{C_{1-\gamma,\rho}[0,h]}\leq
\frac{ h^{\rho(\alpha-\gamma-\lambda+1)}\Gamma(1-\lambda)}{\rho^{\alpha-\gamma+1}\Gamma(\alpha-\lambda+1)}
\big\|\mathcal{M}x_{n}-\mathcal{M}x\big\|_{[0,h]}\label{e3.4}
\end{equation}
Thus, $\|\mathcal{N}x_{n}-\mathcal{N}x\|_{C_{1-\gamma,\rho}[0,h]}\rightarrow0$ as $n\rightarrow+\infty.$ Therefore, $\mathcal{N}$ is continuous.
Moreover, we shall prove that the operator $\mathcal{N}E_{h}$ is equicontinuous. Let $x\in E_{h}$ and $0\leq t_{1}< t_{2}\leq h,$ for any $\delta>0,$ note that
\[\frac{\big(t^{\rho}/\rho\big)^{1-\gamma}\rho^{1-\alpha}}{\Gamma(\alpha)}
\int_{0}^{t}\tau^{\rho-\lambda-1}(t^{\rho}-\tau^{\rho})^{\alpha-1}d\tau=
\frac{ t^{\rho(\alpha-\gamma-\lambda+1)}\Gamma(1-\lambda)}{\rho^{\alpha-\gamma+1}\Gamma(\alpha-\lambda+1)}\rightarrow0\quad as\quad t\rightarrow0^{+},\]
where $0\leq\lambda<1,$ there exists a $(h>)\epsilon_{1}>0$ such that, for $t\in[0,\epsilon_{1}],$ we have
\begin{equation}
\frac{\big(t^{\rho}/\rho\big)^{1-\gamma}\rho^{1-\alpha}L}{\Gamma(\alpha)}
\int_{0}^{t}\tau^{\rho-\lambda-1}(t^{\rho}-\tau^{\rho})^{\alpha-1}d\tau<\frac{\delta}{2}.\label{e3.5}
\end{equation}
In the case, for $t_{1}, t_{2}\in[0,\epsilon_{1}],$ we get
\begin{align}\nonumber
\bigg|\big(t_{1}^{\rho}/\rho\big)^{1-\gamma}&\frac{\rho^{1-\alpha}}{\Gamma(\alpha)}
\int_{0}^{t_{1}}\tau^{\rho-1}(t_{1}^{\rho}-\tau^{\rho})^{\alpha-1}\varphi(\tau,x(\tau))d\tau\\&\nonumber-
\big(t_{2}^{\rho}/\rho\big)^{1-\gamma}\frac{\rho^{1-\alpha}}{\Gamma(\alpha)}
\int_{0}^{t_{2}}\tau^{\rho-1}(t_{2}^{\rho}-\tau^{\rho})^{\alpha-1}\varphi(\tau,x(\tau))d\tau\bigg|\\& \leq\nonumber
\frac{\big(t_{1}^{\rho}/\rho\big)^{1-\gamma}\rho^{1-\alpha}L}{\Gamma(\alpha)}
\int_{0}^{t_{1}}\tau^{\rho-\lambda-1}(t_{_{1}}^{\rho}-\tau^{\rho})^{\alpha-1}d\tau\\&\quad\quad\quad+
\frac{\big(t_{2}^{\rho}/\rho\big)^{1-\gamma}\rho^{1-\alpha}L}{\Gamma(\alpha)}
\int_{0}^{t_{2}}\tau^{\rho-\lambda-1}(t_{2}^{\rho}-\tau^{\rho})^{\alpha-1}d\tau\\&\nonumber<\frac{\delta}{2}+\frac{\delta}{2}=\delta.
\end{align}
In the case, for $t_{1}, t_{2}\in[\frac{\epsilon_{1}}{2},h],$ we have
\begin{align}\nonumber
\bigg|\big(t_{1}^{\rho}/\rho\big)^{1-\gamma}(\mathcal{N}x)(t_{1})&-\big(t_{2}^{\rho}/\rho\big)^{1-\gamma}(\mathcal{N}x)(t_{2})\bigg|
\quad\quad\quad\quad\quad\quad\quad\quad\quad\quad\quad\quad\quad\quad\quad\quad\\&
\nonumber=\bigg|\big(t_{1}^{\rho}/\rho\big)^{1-\gamma}\frac{\rho^{1-\alpha}}{\Gamma(\alpha)}
\int_{0}^{t_{1}}\tau^{\rho-1}(t_{1}^{\rho}-\tau^{\rho})^{\alpha-1}\varphi(\tau,x(\tau))d\tau\\&\nonumber\quad\quad\quad-
\big(t_{2}^{\rho}/\rho\big)^{1-\gamma}\frac{\rho^{1-\alpha}}{\Gamma(\alpha)}
\int_{0}^{t_{2}}\tau^{\rho-1}(t_{2}^{\rho}-\tau^{\rho})^{\alpha-1}\varphi(\tau,x(\tau))d\tau\bigg|\\& \nonumber\leq
\bigg|\frac{\rho^{1-\alpha}}{\Gamma(\alpha)}
\int_{0}^{t_{1}}\tau^{\rho-1}\bigg[\big(t_{1}^{\rho}/\rho\big)^{1-\gamma}(t_{_{1}}^{\rho}-\tau^{\rho})^{\alpha-1}
-\big(t_{2}^{\rho}/\rho\big)^{1-\gamma}(t_{2}^{\rho}-\tau^{\rho})^{\alpha-1}\bigg]\varphi(\tau,x(\tau))d\tau\bigg|\\&\quad\quad\quad+\bigg|
\frac{\big(t_{2}^{\rho}/\rho\big)^{1-\gamma}\rho^{1-\alpha}}{\Gamma(\alpha)}
\int_{t_{1}}^{t_{2}}\tau^{\rho-1}(t_{2}^{\rho}-\tau^{\rho})^{\alpha-1}\varphi(\tau,x(\tau))d\tau\bigg|,
\end{align}
its easy to see form the fact that if $0\leq\nu_{1}<\nu_{2}\leq h,$ then \\
$\big(\nu_{1}^{\rho}/\rho\big)^{1-\gamma}(\nu_{_{1}}^{\rho}-\tau^{\rho})^{\alpha-1}>
\big(\nu_{2}^{\rho}/\rho\big)^{1-\gamma}(\nu_{_{2}}^{\rho}-\tau^{\rho})^{\alpha-1}$ for $0\leq\tau<\nu_{1},$ we get
\begin{align}\nonumber
\bigg|\frac{\rho^{1-\alpha}}{\Gamma(\alpha)}&
\int_{0}^{t_{1}}\tau^{\rho-1}\bigg[\big(t_{1}^{\rho}/\rho\big)^{1-\gamma}(t_{_{1}}^{\rho}-\tau^{\rho})^{\alpha-1}
-\big(t_{2}^{\rho}/\rho\big)^{1-\gamma}(t_{2}^{\rho}-\tau^{\rho})^{\alpha-1}\bigg]\varphi(\tau,x(\tau))d\tau\bigg|\\& \nonumber\leq
\frac{\rho^{1-\alpha}L}{\Gamma(\alpha)}
\int_{0}^{t_{1}}\bigg|\tau^{\rho-\lambda-1}\bigg[\big(t_{1}^{\rho}/\rho\big)^{1-\gamma}(t_{_{1}}^{\rho}-\tau^{\rho})^{\alpha-1}
-\big(t_{2}^{\rho}/\rho\big)^{1-\gamma}(t_{2}^{\rho}-\tau^{\rho})^{\alpha-1}\bigg]\bigg|d\tau\\& \nonumber\leq
\frac{\rho^{1-\alpha}L}{\Gamma(\alpha)}
\int_{0}^{\frac{\epsilon_{1}}{2}}\bigg|\tau^{\rho-\lambda-1}\bigg[\big(t_{1}^{\rho}/\rho\big)^{1-\gamma}(t_{_{1}}^{\rho}-\tau^{\rho})^{\alpha-1}
-\big(t_{2}^{\rho}/\rho\big)^{1-\gamma}(t_{2}^{\rho}-\tau^{\rho})^{\alpha-1}\bigg]\bigg|d\tau\\& \nonumber \quad\quad\quad+
\frac{(\frac{\epsilon_{1}}{2})^{-\lambda}\rho^{1-\alpha}L}{\Gamma(\alpha)}
\int_{\frac{\epsilon_{1}}{2}}^{t_{1}}\bigg|\tau^{\rho-1}\bigg[\big(t_{1}^{\rho}/\rho\big)^{1-\gamma}(t_{_{1}}^{\rho}-\tau^{\rho})^{\alpha-1}
-\big(t_{2}^{\rho}/\rho\big)^{1-\gamma}(t_{2}^{\rho}-\tau^{\rho})^{\alpha-1}\bigg]\bigg|d\tau\\& \nonumber\leq
\frac{2\big((\frac{\epsilon_{1}}{2})^{\rho}/\rho\big)^{1-\gamma}\rho^{1-\alpha}L}{\Gamma(\alpha)}
\int_{0}^{\frac{\epsilon_{1}}{2}}\tau^{\rho-\gamma-1}\big((\frac{\epsilon_{1}}{2})^{\rho}-\tau^{\rho}\big)^{\alpha-1}d\tau\\& \nonumber\quad\quad\quad+
\frac{(\frac{\epsilon_{1}}{2})^{-\lambda}L}{\rho^{\alpha-\gamma+1}\Gamma(\alpha+1)}\bigg[t_{2}^{\rho(1-\gamma)}\big(t_{2}^{\rho}-t_{1}^{\rho}\big)^{\alpha}
-t_{2}^{\rho(1-\gamma)}\big(t_{2}^{\rho}-(\frac{\epsilon_{1}}{2})^{\rho}\big)^{\alpha}+
t_{1}^{\rho(1-\gamma)}\big(t_{1}^{\rho}-(\frac{\epsilon_{1}}{2})^{\rho}\big)^{\alpha}\bigg]\\&\leq
\delta+\frac{(\frac{\epsilon_{1}}{2})^{-\lambda}L}{\rho^{\alpha-\gamma+1}\Gamma(\alpha+1)}\bigg[h^{\rho(1-\gamma)}\big(t_{2}^{\rho}-t_{1}^{\rho}\big)^{\alpha}
+\bigg|t_{2}^{\rho(1-\gamma)}\big(t_{2}^{\rho}-(\frac{\epsilon_{1}}{2})^{\rho}\big)^{\alpha}-
t_{1}^{\rho(1-\gamma)}\big(t_{1}^{\rho}-(\frac{\epsilon_{1}}{2})^{\rho}\big)^{\alpha}\bigg|\bigg].
\end{align}
On the other hand,
\begin{align}\nonumber
\bigg|\frac{\big(t_{2}^{\rho}/\rho\big)^{1-\gamma}\rho^{1-\alpha}}{\Gamma(\alpha)}
\int_{t_{1}}^{t_{2}}\tau^{\rho-1}(t_{2}^{\rho}-&\tau^{\rho})^{\alpha-1}\varphi(\tau,x(\tau))d\tau\bigg|
\quad\quad\quad\quad\quad\quad\quad\quad\quad\quad\quad\quad\quad\\&\nonumber\leq
\frac{(\frac{\epsilon_{1}}{2})^{-\lambda}L}{\rho^{\alpha-\gamma}\Gamma(\alpha)}
\int_{t_{1}}^{t_{2}}\tau^{\rho-1}\big(t_{2}^{\rho}-\tau^{\rho}\big)^{\alpha-1}d\tau
\\&\nonumber=
\frac{(\frac{\epsilon_{1}}{2})^{-\lambda}L}{\rho^{\alpha-\gamma+1}\Gamma(\alpha+1)}\bigg[
t_{2}^{\rho(1-\gamma)}\big(t_{2}^{\rho}-t_{1}^{\rho}\big)^{\alpha}\bigg]\\&\leq
\frac{(\frac{\epsilon_{1}}{2})^{-\lambda}L}{\rho^{\alpha-\gamma+1}\Gamma(\alpha+1)}\bigg[
h^{\rho(1-\gamma)}\big(t_{2}^{\rho}-t_{1}^{\rho}\big)^{\alpha}\bigg].
\end{align}
Obviously, there exists a $(\frac{\epsilon_{1}}{2}>)\epsilon>0$ such that, for
$t_{1}, t_{2}\in[\frac{\epsilon_{1}}{2},h],\quad\big|t_{1}-t_{2}\big|<\epsilon$
implies
\begin{equation}
\big|\big(t_{1}^{\rho}/\rho\big)^{1-\gamma}(\mathcal{N}x)(t_{1})-
\big(t_{2}^{\rho}/\rho\big)^{1-\gamma}(\mathcal{N}x)(t_{2})\big|<2\delta.\label{e3.10}
\end{equation}
Finally, it observe from $(3.6)$ and $(3.10)$ that $\big\{\big(t^{\rho}/\rho\big)^{1-\gamma}\mathcal{N}:x\in E_{h}\big\}$ is equicontinuous.
Evidently, $\big\{\big(t^{\rho}/\rho\big)^{1-\gamma}\mathcal{N}:x\in E_{h}\big\}$ is uniformly bounded, due to $\mathcal{N}E_{h}\subset E_{h}.$
Then, by Lemma 2.5, $\mathcal{N}E_{h}$ is precompact. Thus, $\mathcal{N}$ is completely continuous. Therefore, By Lemma 2.2 (Schauder fixed point theorem) and Lemma 2.4, the initial value problems $(1.3)$ has a local solution.$\quad\quad\quad\quad\quad\quad\Box$
\\\textbf{\ Theorem 3.2.} Assume that a condition $(\mathcal{H}_{2})$  is satisfied. Then the initial value problems $(1.4)$ has at least one solution
$x_{\ell}\in C_{1-\gamma,\rho}[0,h],~(\ell=1,2,...,n)$ for some $(T\geq)h>0.$
\\ \textbf{\ Proof.} Let
\begin{equation}
D=\bigg\{x_{\ell}\in C_{1-\gamma,\rho}[0,T]:\bigg\|x_{\ell}-\frac{x_{0}}{\Gamma(\gamma)}\big(t^{\rho}/\rho\big)^{\gamma-1}\bigg\|_{C_{1-\gamma,\rho}[0,T]}
=\sup_{0\leq t\leq T}\bigg|\big(t^{\rho}/\rho\big)^{1-\gamma}x_{\ell}(t)-\frac{x_{0}}{\Gamma(\gamma)}\bigg|\leq k_{\ell}\bigg\},\label{e3.11}
\end{equation}
where $k_{\ell}>0,~(\ell=1,2,...,n)$ are constants. Since the operators $\mathcal{M}_{\ell},~(\ell=1,2,...,n)$ is bounded, there exists a constant $L_{\ell}>0, (\ell=1,2,...,n)$ such that
\[\sup\big\{\big|(\mathcal{M}_{\ell}x_{\ell})(t)\big|:t\in[0,T],~x_{\ell}\in D\big\}\leq L_{\ell},~(\ell=1,2,...,n).\]
Again, let
\begin{equation}
E_{\ell h}=\bigg\{x_{\ell}:x_{\ell}\in C_{1-\gamma,\rho}[0,T],
\sup_{0\leq t\leq T}\bigg|\big(t^{\rho}/\rho\big)^{1-\gamma}x_{\ell}(t)-\frac{x_{0}}{\Gamma(\gamma)}\bigg|\leq k_{\ell},~\ell=1,2,...,n\bigg\},\label{e3.12}
\end{equation}
where \[h=\min\bigg\{\big(\frac{k_{1}\rho^{\alpha-\gamma_{1}+1}\Gamma(\alpha-\lambda_{1}+1)}{L_{1}\Gamma(1-\lambda_{1})}\big)^{\frac{1}{\rho(\alpha-\gamma-\lambda_{1}+1)}},...,
\big(\frac{k_{n}\rho^{\alpha-\gamma_{n}+1}\Gamma(\alpha-\lambda_{n}+1)}{L_{n}\Gamma(1-\lambda_{n})}\big)^{\frac{1}{\rho(\alpha-\gamma-\lambda_{n}+1)}},T\bigg\}.\]
Obviously, $E_{\ell h}\subseteq C_{1-\gamma,\rho}[0,T]$ be a nonempty, bounded closed and convex subset.
\par Note that $h\leq T$, we can regard $E_{\ell h}$ and $C_{1-\gamma,\rho}[0,T]$ as the restrictions of $D$ and $C_{1-\gamma,\rho}[0,T]$, respectively.
Define the operators $\mathcal{N}_{\ell}$ as follows
\begin{equation}
\left\{\begin{matrix}
(\mathcal{N}_{1}x_{1})(t)=\frac{x_{0}}{\Gamma(\gamma)}\big(t^{\rho}/\rho\big)^{\gamma-1}+\frac{\rho^{1-\alpha}}{\Gamma(\alpha)}
\int_{0}^{t}\tau^{\rho-1}(t^{\rho}-\tau^{\rho})^{\alpha-1}\varphi_{1}(\tau,x_{1}(\tau),x_{2}(\tau),...,x_{n}(\tau))d\tau,
\quad\quad\quad\quad\quad\quad\quad\quad\quad\\\\
(\mathcal{N}_{2}x_{2})(t)=\frac{x_{0}}{\Gamma(\gamma)}\big(t^{\rho}/\rho\big)^{\gamma-1}+\frac{\rho^{1-\alpha}}{\Gamma(\alpha)}
\int_{0}^{t}\tau^{\rho-1}(t^{\rho}-\tau^{\rho})^{\alpha-1}\varphi_{2}(\tau,x_{1}(\tau),x_{2}(\tau),...,x_{n}(\tau))d\tau,
\quad\quad\quad\quad\quad\quad\quad\quad\quad\\ \vdots\\
(\mathcal{N}_{n}x_{n})(t)=\frac{x_{0}}{\Gamma(\gamma)}\big(t^{\rho}/\rho\big)^{\gamma-1}+\frac{\rho^{1-\alpha}}{\Gamma(\alpha)}
\int_{0}^{t}\tau^{\rho-1}(t^{\rho}-\tau^{\rho})^{\alpha-1}\varphi_{n}(\tau,x_{1}(\tau),x_{2}(\tau),...,x_{n}(\tau))d\tau,
\quad\quad\quad\quad\quad\quad\quad\quad\quad\end{matrix}\right.
\end{equation}
for $t\in[0,h].$ Observe that from $(\mathcal{H}_{2})$ and Lemma 2.3, we have $$\mathcal{N}_{\ell}(C_{1-\gamma,\rho}[0,h])\subset C_{1-\gamma,\rho}[0,h],~(\ell=1,2,...,n).$$
\par By relation $(3.13),$ for any $x\in C_{1-\gamma,\rho}[0,h],$ we obtain
\begin{align*}
\bigg|\big(t^{\rho}/\rho\big)^{1-\gamma}(\mathcal{N}_{1}x_{1})(t)-\frac{x_{0}}{\Gamma(\gamma)}\bigg|&=
\bigg|\big(t^{\rho}/\rho\big)^{1-\gamma}\frac{\rho^{1-\alpha}}{\Gamma(\alpha)}
\int_{0}^{t}\tau^{\rho-\lambda_{1}-1}(t^{\rho}-\tau^{\rho})^{\alpha-1}\big[\tau^{\lambda_{1}}
\varphi_{1}(\tau,x_{1}(\tau),x_{2}(\tau),...,x_{n}(\tau))\big]d\tau\bigg|\\& \leq
\frac{\big(t^{\rho}/\rho\big)^{1-\gamma}\rho^{1-\alpha}L_{1}}{\Gamma(\alpha)}
\int_{0}^{t}\tau^{\rho-\lambda_{1}-1}(t^{\rho}-\tau^{\rho})^{\alpha-1}d\tau\\&\leq
\frac{L_{1} \Gamma(1-\lambda_{1})}{\rho^{\alpha-\gamma+1}\Gamma(\alpha-\lambda_{1}+1)}t^{\rho(\alpha-\gamma-\lambda_{1}+1)},
\end{align*}
\begin{align*}
&\bigg|\big(t^{\rho}/\rho\big)^{1-\gamma}(\mathcal{N}_{1}x_{1})(t)-\frac{x_{0}}{\Gamma(\gamma)}\bigg|
\leq\frac{L_{1} \Gamma(1-\lambda_{1})}{\rho^{\alpha-\gamma+1}\Gamma(\alpha-\lambda_{1}+1)}h^{\rho(\alpha-\gamma-\lambda_{1}+1)}\leq k_{1},\\&
\bigg|\big(t^{\rho}/\rho\big)^{1-\gamma}(\mathcal{N}_{2}x_{2})(t)-\frac{x_{0}}{\Gamma(\gamma)}\bigg|
\leq\frac{L_{2} \Gamma(1-\lambda_{2})}{\rho^{\alpha-\gamma+1}\Gamma(\alpha-\lambda_{2}+1)}h^{\rho(\alpha-\gamma-\lambda_{2}+1)}\leq k_{2},\\&
\quad\quad\quad\quad\quad\quad\quad\quad\quad\quad\quad\quad\quad\vdots\\&
\bigg|\big(t^{\rho}/\rho\big)^{1-\gamma}(\mathcal{N}_{n}x_{n})(t)-\frac{x_{0}}{\Gamma(\gamma)}\bigg|
\leq\frac{L_{n} \Gamma(1-\lambda_{n})}{\rho^{\alpha-\gamma+1}\Gamma(\alpha-\lambda_{n}+1)}h^{\rho(\alpha-\gamma-\lambda_{1}+1)}\leq k_{n},
\end{align*}

which yields that $\mathcal{N}_{\ell}E_{\ell h}\subset E_{\ell h},~\ell=1,2,...,n.$
\par Next, we will show that  $\mathcal{N}_{\ell}$ are continuous. For that let $x_{m},x_{\ell}\in E_{\ell h},m>n,~\ell=1,2,...,n$ such that $\|x_{m}-x_{\ell}\|_{C_{1-\gamma,\rho}[0,h]}\rightarrow0$ as $m\rightarrow+\infty.$ In the light of a continuity of $\mathcal{M}_{\ell},$ we have $\|\mathcal{M}_{\ell}x_{m}-\mathcal{M}_{\ell}x_{\ell}\|_{C_{1-\gamma,\rho}[0,h]}\rightarrow0$ as $m\rightarrow+\infty.$
\\Now, Noticing that
\begin{align*}
\bigg|\big(t^{\rho}/\rho\big)^{1-\gamma}(\mathcal{N}_{\ell}x_{m})(t)-&\big(t^{\rho}/\rho\big)^{1-\gamma}(\mathcal{N}_{\ell}x_{\ell})(t)\bigg|\\&=
\bigg|\big(t^{\rho}/\rho\big)^{1-\gamma}\frac{\rho^{1-\alpha}}{\Gamma(\alpha)}
\int_{0}^{t}\tau^{\rho-1}(t^{\rho}-\tau^{\rho})^{\alpha-1}\varphi_{\ell}(\tau,x_{m}(\tau))d\tau\\&\quad\quad\quad\quad-
\big(t^{\rho}/\rho\big)^{1-\gamma}\frac{\rho^{1-\alpha}}{\Gamma(\alpha)}
\int_{0}^{t}\tau^{\rho-1}(t^{\rho}-\tau^{\rho})^{\alpha-1}\varphi_{\ell}(\tau,x_{\ell}(\tau))d\tau\bigg|\\& \leq
\frac{\big(t^{\rho}/\rho\big)^{1-\gamma}\rho^{1-\alpha}}{\Gamma(\alpha)}
\int_{0}^{t}\tau^{\rho-\lambda_{\ell}-1}(t^{\rho}-\tau^{\rho})^{\alpha-1}
\bigg|\tau^{\lambda_{\ell}}\big[\varphi_{\ell}(\tau,x_{m}(\tau))-\varphi_{\ell}(\tau,x_{\ell}(\tau))\big]\bigg|d\tau\\&\leq
\frac{\big(t^{\rho}/\rho\big)^{1-\gamma}\rho^{1-\alpha}}{\Gamma(\alpha)}
\int_{0}^{t}\tau^{\rho-\lambda_{\ell}-1}(t^{\rho}-\tau^{\rho})^{\alpha-1}d\tau
\big\|\mathcal{M}_{\ell}x_{m}-\mathcal{M}_{\ell}x_{\ell}\big\|_{[0,h]}.
\end{align*}
 Then, we have
\begin{equation}
\big\|\mathcal{N}_{\ell}x_{m}-\mathcal{N}_{\ell}x_{\ell}\big\|_{C_{1-\gamma,\rho}[0,h]}\leq
\frac{ h^{\rho(\alpha-\gamma-\lambda_{\ell}+1)}\Gamma(1-\lambda)}{\rho^{\alpha-\gamma+1}\Gamma(\alpha-\lambda_{\ell}+1)}
\big\|\mathcal{M}_{\ell}x_{m}-\mathcal{M}_{\ell}x_{\ell}\big\|_{[0,h]}\label{e3.14}
\end{equation}
Thus, $\|\mathcal{N}_{\ell}x_{m}-\mathcal{N}_{\ell}x_{\ell}\|_{C_{1-\gamma,\rho}[0,h]}\rightarrow0$ as $m\rightarrow+\infty.$ Therefore, $\mathcal{N}_{\ell}$ is continuous.
Moreover, we shall prove that the operators $\mathcal{N}_{\ell}E_{\ell h}$ are equicontinuous. Let $x_{\ell}\in E_{\ell h}$ and $0\leq t_{1}< t_{2}\leq h,$ for any $\delta>0,$ note that
\[\frac{\big(t^{\rho}/\rho\big)^{1-\gamma}\rho^{1-\alpha}}{\Gamma(\alpha)}
\int_{0}^{t}\tau^{\rho-\lambda_{\ell}-1}(t^{\rho}-\tau^{\rho})^{\alpha-1}d\tau=
\frac{ t^{\rho(\alpha-\gamma-\lambda_{\ell}+1)}\Gamma(1-\lambda_{\ell})}{\rho^{\alpha-\gamma+1}\Gamma(\alpha-\lambda_{\ell}+1)}\rightarrow0\quad as\quad t\rightarrow0^{+},\]
where $0\leq\lambda_{\ell}<1,$ there exists a $(h>)\epsilon_{\ell}>0$ such that, for $t\in[0,\epsilon_{\ell}],$ we have
\begin{equation}
\frac{\big(t^{\rho}/\rho\big)^{1-\gamma}\rho^{1-\alpha}L_{\ell}}{\Gamma(\alpha)}
\int_{0}^{t}\tau^{\rho-\lambda_{\ell}-1}(t^{\rho}-\tau^{\rho})^{\alpha-1}d\tau<\frac{\delta}{2}.\label{e3.15}
\end{equation}
In the case, for $t_{1}, t_{2}\in[0,\epsilon_{\ell}],$ we get
\begin{align}\nonumber
\bigg|\big(t_{1}^{\rho}/\rho\big)^{1-\gamma}&\frac{\rho^{1-\alpha}}{\Gamma(\alpha)}
\int_{0}^{t_{1}}\tau^{\rho-1}(t_{1}^{\rho}-\tau^{\rho})^{\alpha-1}\varphi_{\ell}(\tau,x_{\ell}(\tau))d\tau\\&\nonumber-
\big(t_{2}^{\rho}/\rho\big)^{1-\gamma}\frac{\rho^{1-\alpha}}{\Gamma(\alpha)}
\int_{0}^{t_{2}}\tau^{\rho-1}(t_{2}^{\rho}-\tau^{\rho})^{\alpha-1}\varphi_{\ell}(\tau,x_{\ell}(\tau))d\tau\bigg|\\& \leq\nonumber
\frac{\big(t_{1}^{\rho}/\rho\big)^{1-\gamma}\rho^{1-\alpha}L_{\ell}}{\Gamma(\alpha)}
\int_{0}^{t_{1}}\tau^{\rho-\lambda_{\ell}-1}(t_{_{1}}^{\rho}-\tau^{\rho})^{\alpha-1}d\tau\\&\quad\quad\quad+
\frac{\big(t_{2}^{\rho}/\rho\big)^{1-\gamma}\rho^{1-\alpha}L_{\ell}}{\Gamma(\alpha)}
\int_{0}^{t_{2}}\tau^{\rho-\lambda_{\ell}-1}(t_{2}^{\rho}-\tau^{\rho})^{\alpha-1}d\tau\\&\nonumber<\frac{\delta}{2}+\frac{\delta}{2}=\delta.
\end{align}
In the case, for $t_{1}, t_{2}\in[\frac{\epsilon_{\ell}}{2},h],$ we have
\begin{align}\nonumber
\bigg|\big(t_{1}^{\rho}/\rho\big)^{1-\gamma}(\mathcal{N}_{\ell}x_{\ell})(t_{1})&-\big(t_{2}^{\rho}/\rho\big)^{1-\gamma}(\mathcal{N}_{\ell}x_{\ell})(t_{2})\bigg|
\quad\quad\quad\quad\quad\quad\quad\quad\quad\quad\quad\quad\quad\quad\quad\quad\\&
\nonumber=\bigg|\big(t_{1}^{\rho}/\rho\big)^{1-\gamma}\frac{\rho^{1-\alpha}}{\Gamma(\alpha)}
\int_{0}^{t_{1}}\tau^{\rho-1}(t_{1}^{\rho}-\tau^{\rho})^{\alpha-1}\varphi_{\ell}(\tau,x_{\ell}(\tau))d\tau\\&\nonumber\quad\quad\quad-
\big(t_{2}^{\rho}/\rho\big)^{1-\gamma}\frac{\rho^{1-\alpha}}{\Gamma(\alpha)}
\int_{0}^{t_{2}}\tau^{\rho-1}(t_{2}^{\rho}-\tau^{\rho})^{\alpha-1}\varphi_{\ell}(\tau,x_{\ell}(\tau))d\tau\bigg|\\& \nonumber\leq
\bigg|\frac{\rho^{1-\alpha}}{\Gamma(\alpha)}
\int_{0}^{t_{1}}\tau^{\rho-1}\bigg[\big(t_{1}^{\rho}/\rho\big)^{1-\gamma}(t_{_{1}}^{\rho}-\tau^{\rho})^{\alpha-1}
-\big(t_{2}^{\rho}/\rho\big)^{1-\gamma}(t_{2}^{\rho}-\tau^{\rho})^{\alpha-1}\bigg]\varphi_{\ell}(\tau,x_{\ell}(\tau))d\tau\bigg|\\&\quad\quad\quad+\bigg|
\frac{\big(t_{2}^{\rho}/\rho\big)^{1-\gamma}\rho^{1-\alpha}}{\Gamma(\alpha)}
\int_{t_{1}}^{t_{2}}\tau^{\rho-1}(t_{2}^{\rho}-\tau^{\rho})^{\alpha-1}\varphi_{\ell}(\tau,x_{\ell}(\tau))d\tau\bigg|,
\end{align}
its easy to see form the fact that if $0\leq\nu_{1}<\nu_{2}\leq h,$ then \\
$\big(\nu_{1}^{\rho}/\rho\big)^{1-\gamma}(\nu_{_{1}}^{\rho}-\tau^{\rho})^{\alpha-1}>
\big(\nu_{2}^{\rho}/\rho\big)^{1-\gamma}(\nu_{_{2}}^{\rho}-\tau^{\rho})^{\alpha-1}$ for $0\leq\tau<\nu_{1},$ we get
\begin{align}\nonumber
\bigg|\frac{\rho^{1-\alpha}}{\Gamma(\alpha)}&
\int_{0}^{t_{1}}\tau^{\rho-1}\bigg[\big(t_{1}^{\rho}/\rho\big)^{1-\gamma}(t_{_{1}}^{\rho}-\tau^{\rho})^{\alpha-1}
-\big(t_{2}^{\rho}/\rho\big)^{1-\gamma}(t_{2}^{\rho}-\tau^{\rho})^{\alpha-1}\bigg]\varphi_{\ell}(\tau,x_{\ell}(\tau))d\tau\bigg|\\& \nonumber\leq
\frac{\rho^{1-\alpha}L_{\ell}}{\Gamma(\alpha)}
\int_{0}^{t_{1}}\bigg|\tau^{\rho-\lambda_{\ell}-1}\bigg[\big(t_{1}^{\rho}/\rho\big)^{1-\gamma}(t_{_{1}}^{\rho}-\tau^{\rho})^{\alpha-1}
-\big(t_{2}^{\rho}/\rho\big)^{1-\gamma}(t_{2}^{\rho}-\tau^{\rho})^{\alpha-1}\bigg]\bigg|d\tau\\& \nonumber\leq
\frac{\rho^{1-\alpha}L_{\ell}}{\Gamma(\alpha)}
\int_{0}^{\frac{\epsilon_{\ell}}{2}}\bigg|\tau^{\rho-\lambda_{\ell}-1}\bigg[\big(t_{1}^{\rho}/\rho\big)^{1-\gamma}(t_{_{1}}^{\rho}-\tau^{\rho})^{\alpha-1}
-\big(t_{2}^{\rho}/\rho\big)^{1-\gamma}(t_{2}^{\rho}-\tau^{\rho})^{\alpha-1}\bigg]\bigg|d\tau\\& \nonumber \quad\quad\quad+
\frac{(\frac{\epsilon_{\ell}}{2})^{-\lambda_{\ell}}\rho^{1-\alpha}L_{\ell}}{\Gamma(\alpha)}
\int_{\frac{\epsilon_{\ell}}{2}}^{t_{1}}\bigg|\tau^{\rho-1}\bigg[\big(t_{1}^{\rho}/\rho\big)^{1-\gamma}(t_{_{1}}^{\rho}-\tau^{\rho})^{\alpha-1}
-\big(t_{2}^{\rho}/\rho\big)^{1-\gamma}(t_{2}^{\rho}-\tau^{\rho})^{\alpha-1}\bigg]\bigg|d\tau\\& \nonumber\leq
\frac{2\big((\frac{\epsilon_{\ell}}{2})^{\rho}/\rho\big)^{1-\gamma}\rho^{1-\alpha}L_{\ell}}{\Gamma(\alpha)}
\int_{0}^{\frac{\epsilon_{\ell}}{2}}\tau^{\rho-\gamma-1}\big((\frac{\epsilon_{\ell}}{2})^{\rho}-\tau^{\rho}\big)^{\alpha-1}d\tau\\& \nonumber\quad\quad\quad+
\frac{(\frac{\epsilon_{\ell}}{2})^{-\lambda_{\ell}}L_{\ell}}{\rho^{\alpha-\gamma+1}
\Gamma(\alpha+1)}\bigg[t_{2}^{\rho(1-\gamma)}\big(t_{2}^{\rho}-t_{1}^{\rho}\big)^{\alpha}
-t_{2}^{\rho(1-\gamma)}\big(t_{2}^{\rho}-(\frac{\epsilon_{\ell}}{2})^{\rho}\big)^{\alpha}+
t_{1}^{\rho(1-\gamma)}\big(t_{1}^{\rho}-(\frac{\epsilon_{\ell}}{2})^{\rho}\big)^{\alpha}\bigg]\\&\leq
\delta+\frac{(\frac{\epsilon_{\ell}}{2})^{-\lambda}L_{\ell}}{\rho^{\alpha-\gamma+1}\Gamma(\alpha+1)}
\bigg[h^{\rho(1-\gamma)}\big(t_{2}^{\rho}-t_{1}^{\rho}\big)^{\alpha}
+\bigg|t_{2}^{\rho(1-\gamma)}\big(t_{2}^{\rho}-(\frac{\epsilon_{\ell}}{2})^{\rho}\big)^{\alpha}-
t_{1}^{\rho(1-\gamma)}\big(t_{1}^{\rho}-(\frac{\epsilon_{\ell}}{2})^{\rho}\big)^{\alpha}\bigg|\bigg].
\end{align}
On the other hand,
\begin{align}\nonumber
\bigg|\frac{\big(t_{2}^{\rho}/\rho\big)^{1-\gamma}\rho^{1-\alpha}}{\Gamma(\alpha)}
\int_{t_{1}}^{t_{2}}\tau^{\rho-1}(t_{2}^{\rho}-&\tau^{\rho})^{\alpha-1}\varphi_{\ell}(\tau,x_{\ell}(\tau))d\tau\bigg|
\quad\quad\quad\quad\quad\quad\quad\quad\quad\quad\quad\quad\quad\\&\nonumber\leq
\frac{(\frac{\epsilon_{\ell}}{2})^{-\lambda_{\ell}}L_{\ell}}{\rho^{\alpha-\gamma}\Gamma(\alpha)}
\int_{t_{1}}^{t_{2}}\tau^{\rho-1}\big(t_{2}^{\rho}-\tau^{\rho}\big)^{\alpha-1}d\tau
\\&\nonumber=
\frac{(\frac{\epsilon_{\ell}}{2})^{-\lambda_{\ell}}L_{\ell}}{\rho^{\alpha-\gamma+1}\Gamma(\alpha+1)}\bigg[
t_{2}^{\rho(1-\gamma)}\big(t_{2}^{\rho}-t_{1}^{\rho}\big)^{\alpha}\bigg]\\&\leq
\frac{(\frac{\epsilon_{\ell}}{2})^{-\lambda_{\ell}}L_{\ell}}{\rho^{\alpha-\gamma+1}\Gamma(\alpha+1)}\bigg[
h^{\rho(1-\gamma)}\big(t_{2}^{\rho}-t_{1}^{\rho}\big)^{\alpha}\bigg].
\end{align}
Obviously, there exists a $\sigma,~(\frac{\epsilon_{\ell}}{2}>)\sigma>0$ such that, for
$t_{1}, t_{2}\in[\frac{\epsilon_{\ell}}{2},h],\quad\big|t_{1}-t_{2}\big|<\sigma$
implies
\begin{equation}
\big|\big(t_{1}^{\rho}/\rho\big)^{1-\gamma}(\mathcal{N}_{\ell}x_{\ell})(t_{1})-
\big(t_{2}^{\rho}/\rho\big)^{1-\gamma}(\mathcal{N}_{\ell}x_{\ell})(t_{2})\big|<2\delta.\label{e3.20}
\end{equation}
Finally, it observe from $(3.16)$ and $(3.20)$ that $\big\{\big(t^{\rho}/\rho\big)^{1-\gamma}\mathcal{N}_{\ell}:x_{\ell}\in E_{\ell h}\big\}$ is equicontinuous.
Evidently, $\big\{\big(t^{\rho}/\rho\big)^{1-\gamma}\mathcal{N}_{\ell}:x_{\ell}\in E_{\ell h}\big\}$ is uniformly bounded, due to
$\mathcal{N}_{\ell}E_{\ell h}\subset E_{\ell h}.$
Then, by Lemma 2.5, $\mathcal{N}_{\ell}E_{\ell h}$ is precompact. Thus, $\mathcal{N}_{\ell}$ is completely continuous. Therefore, By Lemma 2.2 (Schauder fixed point theorem) and Lemma 2.4, the initial value problems $(1.4)$ has a local solution.$\quad\quad\quad\quad\quad\quad\Box$
\\ \textbf{\ Example 3.1.}
consider the initial value problems
\begin{equation}
\left\{\begin{matrix} _{\frac{1}{2}}D^{\frac{1}{2},\frac{1}{3}}_{0+}x(t)=\varphi(t,x(t)),\quad\quad t\in(0,+\infty),\quad \\\\_{\frac{1}{2}}I^{\frac{1}{3}}_{0+}x(t)\big|_{t=0}=\frac{\sqrt{\pi}}{2},\quad\quad\quad
\quad\quad\quad\quad\quad\quad\quad\end{matrix}\right.
\end{equation}
here,
\par $\alpha=\frac{1}{2},~\beta=\frac{1}{3},~\gamma=\frac{2}{3},~\rho=\frac{1}{2}$ and
$\varphi(t,x(t))=\frac{\sin(1+\sqrt[3]{t^{2}}~x^{2}(t))}{\sqrt[3]{t}},\quad t\in(0,+\infty).$
\\Easily, we can verify that the operator
$(\mathcal{M}x)(t) = t^{\lambda}\varphi(t,x(t))=\sin(1+\sqrt[3]{t^{2}}~x^{2}(t)),$ (where $\lambda=\frac{1}{3},$ in special case),
be the continuous bounded map from $C_{\frac{1}{3},\frac{1}{2}}[0,T]$ into $C[0,T],$ where $T$ is a positive constant.
Then, by Theorem 3.1, the initial value problems $(3.21)$ has a local solution.

\section{The Continuation and Global Existence}

This section contains tow parts, in the first part we discuss the continuation of solution for the initial value problems $(1.3)$
and in the second part we present some results of the global existence. Firstly, we present the following definition and lemma
\\ \textbf{\ Definition 4.1. [20]} Assume that $x(t)~and~\hat{x}(t)$ are solutions of the initial value problems $(1.3)~on~(0,\mu)~and~(0,\hat{\mu}),$
respectively. If $\mu<\hat{\mu}$ and $x(t)=\hat{x}(t)$ for $t\in(0,\mu),$ then we say $\hat{x}(t)$ is the continuation of $x(t),$ or
$x(t)$ can be continued to $(0,\hat{\mu}).$ The solution $x(t)$ is non-continuable if, it has no continuation. The existing
interval of non-continuable solution $x(t)$ is called a maximum existing interval of $x(t).$
\\ \textbf{\ Lemma 4.1. [20]} Let $\rho>0,~h,\mu>0,~0<\alpha<1,~0\leq\nu<1,~\phi_{1}\in C_{\nu}[0,\frac{\mu}{2}]~and~\phi_{2}\in C[\frac{\mu}{2},\mu].$ Then, we have
\[\mathcal{I}_{1}=\int_{0}^{\frac{\mu}{2}}\tau^{\rho-1}(t^{\rho}-\tau^{\rho})^{\alpha-1}\phi_{1}(\tau)d\tau\quad\&\quad
\mathcal{I}_{2}=\int_{\frac{\mu}{2}}^{\mu}\tau^{\rho-1}(t^{\rho}-\tau^{\rho})^{\alpha-1}\phi_{2}(\tau)d\tau \]
are continuous on $[\mu,\mu+h].$
\par Now, we present the first theorem of continuation as follows
\\ \textbf{\ Theorem 4.1.}  Suppose that $(\mathcal{H}_{1})$ is satisfied.
Then, $x=x(t),~t\in(0,\mu)$ is non-continuable if, and only if,
for some $\zeta\in(0,\frac{\mu}{2})$ and any bounded closed subset
$E\subset[\zeta,+\infty)\times\mathbb{R},$ there exists a $t^{*}\in[\zeta,\mu)$
such that, $(t^{*},x(t^{*}))\not\in E.$
\\ \textbf{\ Proof.} Firstly, Assume that $x=x(t)$ be a continuable. Then, there exists solution $\hat{x}(t)$ of initial value problems $(1.3)$
defined on $(0,\hat{\mu})$ such that, $x(t)=\hat{x}(t)$ for $t\in(0,\mu),$ which yields that $\lim_{t\rightarrow\mu^{-}}x(t)=\hat{x}(\mu).$ Now, define
$x(\mu)=\hat{x}(\mu).$ Obviously, $D=\big\{(t,x(t)):t\in[\zeta,\mu)\big\}$ be a compact subset of $[\zeta,+\infty)\times\mathbb{R}.$ Moreover, there exists on
$t^{*}\in[\zeta,\mu)$ such that, $(t^{*},x(t^{*}))\not\in D.$ This contradiction gives that $x(t)$ is non-continuable.
\par Secondly, Assume that there exists a compact
subset $E\subset[\zeta,+\infty)\times\mathbb{R}$ such that, $\big\{(t,x(t)):t\in[\zeta,\mu)\big\}\subset E.$ Then, a compactness of $E$ yields that
$\mu<+\infty.$ By $(\mathcal{H}_{1}),$ there exists $\Lambda>0$ such that,
$\sup_{(t,x(t))\in E}\big|\varphi(t,x(t))\big|\leq\Lambda$
\\\textbf{Step: 1.} Now, we show that $\lim_{t\rightarrow\mu^{-}}x(t)$ exists. For that we put
\begin{align}
&\Psi(\tau,t)=\bigg|\frac{x_{0}}{\Gamma(\gamma)}\big(\tau^{\rho}/\rho\big)^{\gamma-1}-
\frac{x_{0}}{\Gamma(\gamma)}\big(t^{\rho}/\rho\big)^{\gamma-1}\bigg|,\quad(\tau,t)\in[2\zeta,\mu]\times[2\zeta,\mu]\\&
\mathcal{J}=\int_{0}^{\zeta}\tau^{\rho-\lambda-1}(t^{\rho}-\tau^{\rho})^{\alpha-1}(\tau)d\tau,\quad t\in[2\zeta,\mu].
\end{align}
Easily, we can verify that $\Psi(\tau,t)~and~\mathcal{J}$ are uniformly continuous on $[2\zeta,\mu]\times[2\zeta,\mu]~and~[2\zeta,\mu],$ respectively.
\par Next, $\forall t_{1},t_{2}\in[2\zeta,\mu],~t_{1}<t_{2},$ by using equation $(4.1)$ we have
\begin{align}\nonumber
\big|x(t_{1})-x(t_{2})\big|&\\&\nonumber=\bigg|\frac{x_{0}}{\Gamma(\gamma)}\big(t_{1}^{\rho}/\rho\big)^{\gamma-1}+\frac{\rho^{1-\alpha}}{\Gamma(\alpha)}
\int_{0}^{t_{1}}\tau^{\rho-1}(t_{1}^{\rho}-\tau^{\rho})^{\alpha-1}\varphi(\tau,x(\tau))d\tau\\&\nonumber\quad\quad\quad\quad\quad\quad\quad
-\bigg[\frac{x_{0}}{\Gamma(\gamma)}\big(t_{2}^{\rho}/\rho\big)^{\gamma-1}+\frac{\rho^{1-\alpha}}{\Gamma(\alpha)}
\int_{0}^{t_{2}}\tau^{\rho-1}(t_{2}^{\rho}-\tau^{\rho})^{\alpha-1}\varphi(\tau,x(\tau))d\tau\bigg]\bigg|
\\&
\nonumber\leq\Psi(t_{1},t_{2})+\bigg|\frac{\rho^{1-\alpha}}{\Gamma(\alpha)}
\int_{0}^{t_{1}}\tau^{\rho-1}(t_{1}^{\rho}-\tau^{\rho})^{\alpha-1}\varphi(\tau,x(\tau))d\tau\\&\nonumber\quad\quad\quad\quad\quad\quad\quad-
\frac{\rho^{1-\alpha}}{\Gamma(\alpha)}
\int_{0}^{t_{2}}\tau^{\rho-1}(t_{2}^{\rho}-\tau^{\rho})^{\alpha-1}\varphi(\tau,x(\tau))d\tau\bigg|\\& \nonumber\leq
\Psi(t_{1},t_{2})+\bigg|\frac{\rho^{1-\alpha}}{\Gamma(\alpha)}
\int_{0}^{\zeta}\tau^{\rho-\lambda-1}\bigg[(t_{_{1}}^{\rho}-\tau^{\rho})^{\alpha-1}
-(t_{2}^{\rho}-\tau^{\rho})^{\alpha-1}\bigg](\mathcal{M}x)(\tau)d\tau\bigg|\\&\nonumber\quad\quad\quad\quad\quad\quad+
\frac{\rho^{1-\alpha}}{\Gamma(\alpha)}
\int_{\zeta}^{t_{1}}\tau^{\rho-1}\bigg[(t_{_{1}}^{\rho}-\tau^{\rho})^{\alpha-1}
-(t_{2}^{\rho}-\tau^{\rho})^{\alpha-1}\bigg]\big|\varphi(\tau,x(\tau))\big|d\tau\\&\nonumber\quad\quad\quad\quad\quad\quad+
\frac{\rho^{1-\alpha}}{\Gamma(\alpha)}
\int_{t_{1}}^{t_{2}}\tau^{\rho-1}(t_{2}^{\rho}-\tau^{\rho})^{\alpha-1}\big|\varphi(\tau,x(\tau))\big|d\tau,
\\& \nonumber\leq
\Psi(t_{1},t_{2})+\frac{\|\mathcal{M}x\|_{[0,\zeta]}\rho^{1-\alpha}}{\Gamma(\alpha)}
\int_{0}^{\zeta}\tau^{\rho-\lambda-1}\bigg|\bigg[(t_{_{1}}^{\rho}-\tau^{\rho})^{\alpha-1}
-(t_{2}^{\rho}-\tau^{\rho})^{\alpha-1}\bigg]\bigg|d\tau\\&\nonumber\quad\quad\quad\quad\quad\quad+
\frac{\Lambda\rho^{1-\alpha}}{\Gamma(\alpha)}
\int_{\zeta}^{t_{1}}\tau^{\rho-1}\bigg[(t_{_{1}}^{\rho}-\tau^{\rho})^{\alpha-1}
-(t_{2}^{\rho}-\tau^{\rho})^{\alpha-1}\bigg]d\tau\\&\nonumber\quad\quad\quad\quad\quad\quad+
\frac{\Lambda\rho^{1-\alpha}}{\Gamma(\alpha)}
\int_{t_{1}}^{t_{2}}\tau^{\rho-1}(t_{2}^{\rho}-\tau^{\rho})^{\alpha-1}d\tau\\&\nonumber\leq
\Psi(t_{1},t_{2})+\frac{\|\mathcal{M}x\|_{[0,\zeta]}\rho^{1-\alpha}}{\Gamma(\alpha)}
\big|\mathcal{J}(t_{1})-\mathcal{J}(t_{2})\big|\\&\quad\quad\quad\quad\quad\quad+
\frac{\Lambda}{\rho^{\alpha}\Gamma(\alpha+1)}\bigg[2(t_{2}^{\rho}-t_{1}^{\rho})^{\alpha}+(t_{1}^{\rho}-\zeta^{\rho})^{\alpha}-
(t_{2}^{\rho}-\zeta^{\rho})^{\alpha}\bigg].
\end{align}
By uniform continuity of $\Psi(\tau,t)~and~\mathcal{J}(t),$ together with a Cauchy convergence criterion, we get
$\lim_{t\rightarrow\mu^{-}}x(t)=x^{*}.$
\\\textbf{Step: 2.} In this Part we show that $x(t)$ is a continuable. Since $E$ be the closed subset, we have $(\mu,x^{*})\in E.$ Define
$x(\mu)=x^{*},$ then $x(t)\in C_{1-\gamma,\rho}[0,\mu].$ We denote
\begin{equation}
x_{1}(t)=\frac{x_{0}}{\Gamma(\gamma)}\big(t^{\rho}/\rho\big)^{\gamma-1}+\frac{\rho^{1-\alpha}}{\Gamma(\alpha)}
\int_{0}^{\mu}\tau^{\rho-1}(t^{\rho}-\tau^{\rho})^{\alpha-1}\varphi(\tau,x(\tau))d\tau,\quad t\in[\mu,\mu+1],
\end{equation}
and we define the operator $\mathcal{K}$ as follows
\begin{equation}
(\mathcal{K}y)(t)=x_{1}(t)+\frac{\rho^{1-\alpha}}{\Gamma(\alpha)}
\int_{\mu}^{t}\tau^{\rho-1}(t^{\rho}-\tau^{\rho})^{\alpha-1}\varphi(\tau,y(\tau))d\tau,\quad t\in[\mu,\mu+1],
\end{equation}
where $y\in C[\mu,\mu+1].$ In the light of Lemmas 2.3 and 4.1, we get $\mathcal{K}(C[\mu,\mu+1])\subset C[\mu,\mu+1].$
\par Now, assume that
\begin{equation}
S_{k}=\bigg\{(t,y): \mu\leq t\leq\mu+1,~|y|\leq\max_{\mu\leq t\leq\mu+1}|x(t)|+k\bigg\},\quad k>0.
\end{equation}
In the view of a continuity of $\varphi~on~S_{k},$ we can denote $\Theta=\max_{(t,y)\in S_{k}}|\varphi(t,y)|.$
\par Again, assume that
\begin{equation}
S_{h}=\bigg\{y\in[\mu,\mu+h]: \max_{t [\mu,\mu+h]}|y(t)-x_{1}(t)|\leq k,\quad y(\mu)=x_{1}(\mu)\bigg\},
\end{equation}
where $h=\min\bigg\{\bigg(\frac{k\rho^{\alpha}\Gamma(\alpha+1)}{\Theta}\bigg)^{\frac{1}{\rho\alpha}},1\bigg\}.$ We can claim
that the operator $\mathcal{K}$ is a completely continuous on $S_{h}.$ Firstly, we will show that  $\mathcal{K}$ is a continuous. Put
$\{y_{n}\}\subseteq C[\mu,\mu+h],~\|y_{n}-y\|\rightarrow0~as~n\rightarrow+\infty.$ So, we have
\begin{align}\nonumber
\big|(\mathcal{K}y_{n})(t)-(\mathcal{K}y)(t)\big|&\\&\nonumber=
\bigg|\frac{\rho^{1-\alpha}}{\Gamma(\alpha)}\int_{\mu}^{t}\tau^{\rho-1}(t^{\rho}-\tau^{\rho})^{\alpha-1}
\big[\varphi(\tau,y_{n}(\tau))-\varphi(\tau,y(\tau))\big]d\tau\bigg|
\\&\nonumber\leq\big\|\varphi(\tau,y_{n}(\tau))-\varphi(\tau,y(\tau))\big\|_{[\mu,\mu+h]}
\frac{\rho^{1-\alpha}}{\Gamma(\alpha)}\int_{\mu}^{t}\tau^{\rho-1}(t^{\rho}-\tau^{\rho})^{\alpha-1}d\tau
\\&\leq\big\|\varphi(\tau,y_{n}(\tau))-\varphi(\tau,y(\tau))\big\|_{[\mu,\mu+h]}
\frac{h^{\rho\alpha}}{\rho^{\alpha}\Gamma(\alpha+1)}.
\end{align}
By a continuity of $\varphi~on~S_{k},$ we get $\big\|\varphi(\tau,y_{n}(\tau))-\varphi(\tau,y(\tau))\big\|_{[\mu,\mu+h]}
\rightarrow0~as~n\rightarrow+\infty.$ Thus, $\big\|\mathcal{K} y_{n}-\mathcal{K}y\big\|_{[\mu,\mu+h]}
\rightarrow0~as~n\rightarrow+\infty,$ which yields that the $\mathcal{K}$ is a continuous.
\par Secondly, we will prove that $\mathcal{K}S_{h}$ is equicontinuous. For any $y\in S_{h},$ we have $(\mathcal{K}y)(\mu)=x_{1}(\mu)$ and
\begin{align}\nonumber
\big|(\mathcal{K}y)(t)-x_{1}(t)\big|&=
\bigg|\frac{\rho^{1-\alpha}}{\Gamma(\alpha)}\int_{\mu}^{t}\tau^{\rho-1}(t^{\rho}-\tau^{\rho})^{\alpha-1}
\varphi(\tau,y(\tau))d\tau\bigg|
\\&\nonumber\leq\frac{\rho^{1-\alpha}}{\Gamma(\alpha)}\int_{\mu}^{t}\tau^{\rho-1}(t^{\rho}-\tau^{\rho})^{\alpha-1}\big|\varphi(\tau,y(\tau))\big|d\tau
\\&\leq\frac{\Theta (t^{\rho}-\mu^{\rho})^{\alpha}}{\rho^{\alpha}\Gamma(\alpha+1)}\leq\frac{h^{\rho\alpha}}{\rho^{\alpha}\Gamma(\alpha+1)}\leq k.
\end{align}
Therefore, $\mathcal{K}S_{h}\subset S_{h}.$
\par Now, put
\[\mathcal{I}(t)=\frac{\rho^{1-\alpha}}{\Gamma(\alpha)}\int_{0}^{\mu}\tau^{\rho-1}(t^{\rho}-\tau^{\rho})^{\alpha-1}
\varphi(\tau,x(\tau))d\tau.\]
By using Lemma 4.1, $\mathcal{I}(t)$ is a continuous on $[\mu,\mu+h]. \forall y\in S_{h},~\mu\leq t_{1}<t_{2}\leq\mu+h,$ we have
\begin{align}\nonumber
\big|(\mathcal{K}y)(t_{1})-(\mathcal{K}y)(t_{2})\big|&\\&\nonumber\leq
\Psi(t_{1},t_{2})+\bigg|\frac{\rho^{1-\alpha}}{\Gamma(\alpha)}
\int_{0}^{\mu}\tau^{\rho-1}\bigg[(t_{_{1}}^{\rho}-\tau^{\rho})^{\alpha-1}
-(t_{2}^{\rho}-\tau^{\rho})^{\alpha-1}\bigg]\varphi(\tau,x(\tau))d\tau\bigg|\\&\nonumber\quad\quad\quad\quad\quad\quad+\bigg|
\frac{\rho^{1-\alpha}}{\Gamma(\alpha)}
\int_{\mu}^{t_{1}}\tau^{\rho-1}\bigg[(t_{_{1}}^{\rho}-\tau^{\rho})^{\alpha-1}
-(t_{2}^{\rho}-\tau^{\rho})^{\alpha-1}\bigg]\varphi(\tau,x(\tau))d\tau\bigg|\\&\nonumber\quad\quad\quad\quad\quad\quad+\bigg|
\frac{\rho^{1-\alpha}}{\Gamma(\alpha)}
\int_{t_{1}}^{t_{2}}\tau^{\rho-1}(t_{2}^{\rho}-\tau^{\rho})^{\alpha-1}\varphi(\tau,x(\tau))d\tau\bigg|,
\\& \nonumber\leq
\Psi(t_{1},t_{2})+\big|\mathcal{I}(t_{1})-\mathcal{I}(t_{2})\big|\\&\nonumber\quad\quad\quad\quad\quad\quad+
\frac{\rho^{1-\alpha}}{\Gamma(\alpha)}
\int_{\mu}^{t_{1}}\tau^{\rho-1}\bigg[(t_{_{1}}^{\rho}-\tau^{\rho})^{\alpha-1}
-(t_{2}^{\rho}-\tau^{\rho})^{\alpha-1}\bigg]\big|\varphi(\tau,x(\tau))\big|d\tau\\&\nonumber\quad\quad\quad\quad\quad\quad+
\frac{\rho^{1-\alpha}}{\Gamma(\alpha)}
\int_{t_{1}}^{t_{2}}\tau^{\rho-1}(t_{2}^{\rho}-\tau^{\rho})^{\alpha-1}\big|\varphi(\tau,x(\tau))\big|d\tau\\&\nonumber\leq
\Psi(t_{1},t_{2})+\big|\mathcal{I}(t_{1})-\mathcal{I}(t_{2})\big|\\&\quad\quad\quad\quad\quad\quad+
\frac{\Lambda}{\rho^{\alpha}\Gamma(\alpha+1)}\bigg[2(t_{_{2}}^{\rho}-t_{1}^{\rho})^{\alpha}+(t_{1}^{\rho}-\mu^{\rho})^{\alpha}-
(t_{2}^{\rho}-\mu^{\rho})^{\alpha}\bigg].
\end{align}
In the view of uniform continuity of $\mathcal{I}(t)~on~[\mu,\mu+h]$ and inequality $(4.10),$ we conclude that $\big\{(\mathcal{K}y)(t):y\in S_{h}\big\}$
is equicontinuous. Thus, the operator $\mathcal{K}$ is completely continuous. Therefore, By Lemma 2.2 (Schauder fixed point theorem)
$\mathcal{K}$ has a fixed point $\hat{x}(t)\in S_{h},~i.e.$
\begin{align}\nonumber
\hat{x}(t)&=x_{1}(t)+\frac{\rho^{1-\alpha}}{\Gamma(\alpha)}
\int_{\mu}^{t}\tau^{\rho-1}(t^{\rho}-\tau^{\rho})^{\alpha-1}\varphi(\tau,\hat{x}(\tau))d\tau\\&=
\frac{x_{0}}{\Gamma(\gamma)}\big(t^{\rho}/\rho\big)^{1-\gamma}+\frac{\rho^{1-\alpha}}{\Gamma(\alpha)}
\int_{0}^{t}\tau^{\rho-1}(t^{\rho}-\tau^{\rho})^{\alpha-1}\varphi(\tau,\tilde{x}(\tau))d\tau,\quad t\in[\mu,\mu+h]
\end{align}
where,
\begin{equation*}
\tilde{x}(t)=\left\{\begin{matrix} x(t)\quad\quad if \quad t\in(0,\mu],~~\quad
\\\\ \hat{x}(t)\quad\quad if \quad t\in[\mu,\mu+h].\end{matrix}\right.
\end{equation*}
From Lemma 2.1, it follows that $\tilde{x}\in C_{1-\gamma,\rho}[o,\mu+h]$ and
\begin{equation*}
\tilde{x}(t)=\frac{x_{0}}{\Gamma(\gamma)}\big(t^{\rho}/\rho\big)^{\gamma-1}+\frac{\rho^{1-\alpha}}{\Gamma(\alpha)}
\int_{0}^{t}\tau^{\rho-1}(t^{\rho}-\tau^{\rho})^{\alpha-1}\varphi(\tau,\tilde{x}(\tau))d\tau.
\end{equation*}
Hence, In the light of Lemma 2.3, $\tilde{x}(t)$ is a solution of the initial value problems $(1.3)~on~[o,\mu+h].$ This gives contradiction because $x(t)$
is non-continuable.$\quad\quad\quad\Box$
\par Now, we will give the second theorem of continuation, which is a more  applied of convenient.
\\ \textbf{\ Theorem 4.2.}  Suppose that $(\mathcal{H}_{1})$ is satisfied.
Then, $x=x(t),~t\in(0,\mu)$ is non-continuable if, and only if,
\begin{equation}
\lim_{t\rightarrow\mu^{-}}\sup|\Phi(t)|=+\infty,
\end{equation}
where $\quad\Phi(t)=(t,x(t))\quad\&\quad|\Phi(t)|=\sqrt{t^{2}+x^{2}(t)}.$
\\ \textbf{\ Proof.} Firstly, Assume that $x=x(t)$ be a continuable. Then, there exists solution $\hat{x}(t)$ of initial value problems $(1.3)$
defined on $(0,\hat{\mu})$ such that, $x(t)=\hat{x}(t)$ for $t\in(0,\mu),$ which yields that $\lim_{t\rightarrow\mu^{-}}x(t)=\hat{x}(\mu).$ Therefore,
$|\Phi(t)|\rightarrow|\Phi(\mu)|~as~t\rightarrow\mu^{-},$ which gives a contradiction.
\par Secondly, Assume that equation $(4.12)$ is not true. Then, there exist a sequence $\{t_{m}\}~and~M>0,$ where $M$ is positive constant, such that
\begin{align}\nonumber
&\quad\quad t_{m}<t_{m+1},~~m\in\mathbb{N},\\&\quad\quad
\lim_{m\rightarrow\infty}t_{m}=\mu,~~|\Phi(t_{m})|\leq M,\\&\nonumber
i.e.\quad t_{m}^{2}+x^{2}(t_{m})\leq M^{2}.
\end{align}
Since $x(t_{m})$ be the bounded convergent subsequence, without loss of generality,
we put
\begin{equation}
\lim_{m\rightarrow\infty}x(t_{m})=x^{*}.
\end{equation}
Now, for any given $\delta>0,$ there exists $T\in(0,\mu)$ such that, $\big|x(t)-x^{*}\big|<\delta,~t\in(T,\mu),$ we show that
\begin{equation}
\lim_{t\rightarrow\mu^{-}}x(t_{m})=x^{*}.
\end{equation}
For sufficiently small $\zeta>0,$ let
\begin{equation}
S_{1}=\bigg\{(t,x): t\in [\zeta,\mu],~~|x|\leq \sup_{t\in [\zeta,\mu)}|x(t)|\bigg\}.
\end{equation}
In the light of continuity of $\varphi~on~S_{1},$ we denote $\Phi=\max_{(t,y)\in S_{1}}|\varphi(t,y)|.$ From equations $(4.13)~and~(4.14),$ it follows that
there exists $m_{0}$ such that $t_{m_{0}}>\zeta$ and for $m\geq m_{0},$ we have
\begin{equation*}
\big|x(t_{m})-x^{*}\big|\leq\frac{\delta}{2}.
\end{equation*}
If $(4.14)$ is not true, then for $m\geq m_{0},$ there exists $\xi_{m}\in(t_{m},\mu)$ such that, for
$t\in(t_{m},\xi_{m}),~\big|x(t)-x^{*}\big|<\delta~and~\big|x(\xi_{m})-x^{*}\big|\geq\delta.$ Hence,
\begin{align}\nonumber
\delta&\leq\big|x(\xi_{m})-x^{*}\big|\leq \big|x(t_{m})-x^{*}\big|+\big|x(\xi_{m})-x(t_{m})\big|\\&\nonumber\leq\frac{\delta}{2}+
\bigg|\frac{\rho^{1-\alpha}}{\Gamma(\alpha)}
\int_{0}^{t_{m}}\tau^{\rho-1}(t_{m}^{\rho}-\tau^{\rho})^{\alpha-1}\varphi(\tau,x(\tau))d\tau-\frac{\rho^{1-\alpha}}{\Gamma(\alpha)}
\int_{0}^{\xi_{m}}\tau^{\rho-1}(\xi_{m}^{\rho}-\tau^{\rho})^{\alpha-1}\varphi(\tau,x(\tau))d\tau\bigg|\\&\nonumber\leq\frac{\delta}{2}+
\frac{\rho^{1-\alpha}}{\Gamma(\alpha)}\bigg|\int_{0}^{\zeta}\tau^{\rho-1}
\bigg[(t_{m}^{\rho}-\tau^{\rho})^{\alpha-1}-(\xi_{m}^{\rho}-\tau^{\rho})^{\alpha-1}\bigg]\varphi(\tau,x(\tau))d\tau\bigg|\\&\nonumber\quad\quad+
\frac{\rho^{1-\alpha}}{\Gamma(\alpha)}\bigg|\int_{\zeta}^{t_{m}}\tau^{\rho-1}
\bigg[(t_{m}^{\rho}-\tau^{\rho})^{\alpha-1}-(\xi_{m}^{\rho}-\tau^{\rho})^{\alpha-1}\bigg]\varphi(\tau,x(\tau))d\tau\bigg|\\&\nonumber\quad\quad+
\frac{\rho^{1-\alpha}}{\Gamma(\alpha)}\bigg|\int_{t_{m}}^{\xi_{m}}\tau^{\rho-1}(\xi_{m}^{\rho}-\tau^{\rho})^{\alpha-1}\varphi(\tau,x(\tau))d\tau\bigg|
\\&\nonumber\leq\frac{\delta}{2}+\frac{\|\mathcal{M}x\|_{[0,\zeta]}\rho^{1-\alpha}}{\Gamma(\alpha)}
\big|\mathcal{J}(t_{m})-\mathcal{J}(\xi_{m})\big|\\&\quad\quad\quad\quad\quad\quad+
\frac{\Phi}{\rho^{\alpha}\Gamma(\alpha+1)}\bigg[2(\xi_{m}^{\rho}-t_{m}^{\rho})^{\alpha}+(t_{m}^{\rho}-\zeta^{\rho})^{\alpha}-
(\xi_{m}^{\rho}-\zeta^{\rho})^{\alpha}\bigg],
\end{align}
where $\mathcal{J}(t)$ is defined  in the equation $(4.2).$ By a continuity of $\mathcal{J}(t)~on~[t_{m_{0}},\mu],$ and for sufficiently large $m\geq m_{0},$
we have
\begin{equation}
\frac{\|\mathcal{M}x\|_{[0,\zeta]}\rho^{1-\alpha}}{\Gamma(\alpha)}
\big|\mathcal{J}(t_{m})-\mathcal{J}(\xi_{m})\big|+
\frac{\Phi}{\rho^{\alpha}\Gamma(\alpha+1)}\bigg[2(\xi_{m}^{\rho}-t_{m}^{\rho})^{\alpha}+(t_{m}^{\rho}-\zeta^{\rho})^{\alpha}-
(\xi_{m}^{\rho}-\zeta^{\rho})^{\alpha}\bigg]<\frac{\delta}{2}.
\end{equation}
From equations $(4.17)~and~(4.18),$ we obtain
\begin{equation*}
\delta\leq\big|x(\xi_{m})-x^{*}\big|<\frac{\delta}{2}+\frac{\delta}{2}=\delta.
\end{equation*}
This contradiction gives that $\lim_{t\rightarrow\mu^{-}}x(t)$ exists.
\par By using same argument as in a proof of the previous theorem, easily we can prove  the continuation of
$x(t).\quad\quad\quad\quad\quad\quad\quad\quad\quad\quad\quad\quad\quad\quad\quad\quad\quad\quad\quad\quad\quad\Box$
\par In the next part of this section, we discuss the global existence of solutions for the initial value problems $(1.3),$
which is based on results obtained previously. Applying the second theorem of continuation (Theorem 4.2),
we can immediately obtain the following conclusion about the global existence of solution for the initial value problems $(1.3).$
\\ \textbf{\ Theorem 4.3.}  Suppose that $(\mathcal{H}_{1})$ is satisfied. Let $x(t)$ is the solution of the initial value problems
$(1.3)~on~(0,\mu).$ If $x(t)$ be a bounded on $[\zeta,\mu)$ for some $\zeta>0,$ then $\mu=+\infty.$
\par For illustrative the above theorem we give the following example
\\ \textbf{\ Example 4.1.} We consider the initial value problem as following
\begin{equation}
\left\{\begin{matrix} _{2}D^{\frac{1}{2},\frac{1}{4}}_{0+}x(t)=\varphi(t,x(t)),\quad\quad t\in(0,+\infty),\quad \\\\_{2}I^{\frac{3}{8}}_{0+}x(t)\big|_{t=0}=1,\quad\quad\quad\quad
\quad\quad\quad\quad\quad\quad\quad\end{matrix}\right.
\end{equation}
here,
\par $\alpha=\frac{1}{2},~\beta=\frac{1}{4},~\gamma=\frac{5}{8},~\rho=2$ and
$\varphi(t,x(t))=\frac{\exp(-t^{2}x\sin t)}{\sqrt{t}~(1-t)}.$
\\By applying Theorem 3.1, we know that the initial value problems $(4.19)$
has at least one a local solution $x(t)~on~(0,h]$ for some $h>0.$ In the view of the Lemma 2.4,
$x(t)$ satisfies the following integral equation
\begin{equation}\quad\quad\quad\quad
x(t)=\frac{1}{\Gamma(\frac{5}{8})}\big(t^{2}/2\big)^{-\frac{3}{8}}+\frac{\sqrt{2}}{\Gamma(\frac{1}{2})}
\int_{0}^{t}\tau~\frac{\exp(-\tau^{2}x\sin\tau)}{\sqrt{\tau(t^{2}-\tau^{2})}~(1-\tau)}d\tau.
\end{equation}
Hence,
\begin{equation}\quad\quad\quad\quad
|x(t)|\leq\frac{1}{\Gamma(\frac{5}{8})}\big(t^{2}/2\big)^{-\frac{3}{8}}+\frac{\sqrt{2}}{\Gamma(\frac{1}{2})}
\int_{0}^{t}\frac{\tau}{\sqrt{\tau(t^{2}-\tau^{2})}}d\tau.
\end{equation}
Assume that $[0,\mu)~with~( \mu<+\infty )$ be a maximum existing interval of $x(t).$ Easily, we can see that for any $\zeta\in (0,\mu), x(t)$ be a
bounded on $[\zeta,\mu).$ By using Theorem 4.3, we have $\mu=+\infty,~i.e.$ a maximum existing interval
of $x(t)$ is $(0,+\infty).$
\par In the light of (Lemma 7.1.1, in [23], Theorem 1, in [24] and Lemma 7.14, in [25]), we state a more generalization of Gronwall's lemma
for singular kernels which is essential for our discussion.
\\ \textbf{\ Lemma 4.2.} Assume that $\phi:[o,\mu]\rightarrow[o,+\infty)$ is a real function and $\theta(.)$ is a
non-negative locally integrable function on $[0,\mu].$ And let there exists $\rho>0,~\omega>0,~and~0<\alpha<1,$ such that
\begin{equation*}\quad\quad\quad\quad
\phi(t)\leq\theta(t)+\omega\int_{0}^{t}\tau^{\rho-1}\frac{\rho^{\alpha}}{(t^{\rho}-\tau^{\rho})^{\alpha}}\phi(\tau)d\tau.
\end{equation*}
Then, there exists a constant $C=C(\alpha)$ such that for $t\in[0,\mu],$ we have
\begin{equation*}\quad\quad\quad\quad
\phi(t)\leq\theta(t)+C~\omega\int_{0}^{t}\tau^{\rho-1}\frac{\rho^{\alpha}}{(t^{\rho}-\tau^{\rho})^{\alpha}}\phi(\tau)d\tau.
\end{equation*}
\\ \textbf{\ Theorem 4.4.} Suppose that $(\mathcal{H}_{1})$ is satisfied and there exist three non-negative continuous functions 
$f(t),~g(t),~\psi(t):[0,+\infty)\rightarrow[0,+\infty)$ such that 
$\big|\varphi(t,x(t))\big|\leq g(t) f(|x(t)|)+\psi(t),~where~g(\eta)\leq\eta~for~\eta\geq0.$ 
Then the initial value problems $(1.3)$ has one solution in $C_{1-\gamma,\rho}[0,+\infty).$
\\ \textbf{\ Proof.} The local existence of a solution of the initial value problems $(1.3)$ can be deduced from Theorem 3.1. 
By using Lemma 2.4, $x(t)$ satisfies the second kind Volterra fractional integral equation
\begin{equation}
x(t)=\frac{x_{0}}{\Gamma(\gamma)}\big(t^{\rho}/\rho\big)^{\gamma-1}+\frac{\rho^{1-\alpha}}{\Gamma(\alpha)}
\int_{0}^{t}\frac{\tau^{\rho-1}\varphi(\tau,x(\tau))}{(t^{\rho}-\tau^{\rho})^{1-\alpha}}d\tau.
\end{equation} 
\par Assume that $[0,\mu)~with~( \mu<+\infty )$ be a maximum existing interval of $x(t).$ Then, we have
\begin{align}\nonumber
\big|\big(t^{\rho}/\rho\big)^{1-\gamma}x(t)\big|&=\bigg|\frac{x_{0}}{\Gamma(\gamma)}+\big(t^{\rho}/\rho\big)^{1-\gamma}\frac{\rho^{1-\alpha}}{\Gamma(\alpha)}
\int_{0}^{t}\tau^{\rho-1}(t^{\rho}-\tau^{\rho})^{\alpha-1}\varphi(\tau,x(\tau))d\tau\bigg|\\&\nonumber\leq
\frac{x_{0}}{\Gamma(\gamma)}+\frac{\big(t^{\rho}/\rho\big)^{1-\gamma}\rho^{1-\alpha}}{\Gamma(\alpha)}
\int_{0}^{t}\tau^{\rho-1}(t^{\rho}-\tau^{\rho})^{\alpha-1}\big|\varphi(\tau,x(\tau))\big|d\tau\\&\nonumber\leq
\frac{x_{0}}{\Gamma(\gamma)}+\frac{\mu^{\rho(1-\gamma)}}{\rho^{\alpha-\gamma}\Gamma(\alpha)}
\int_{0}^{t}\tau^{\rho-1}(t^{\rho}-\tau^{\rho})^{\alpha-1}\big[g(\tau) f((\tau^{\rho}/\rho)^{1-\gamma}|x(\tau)|)+\psi(\tau)\big]d\tau\\&\nonumber\leq
\frac{x_{0}}{\Gamma(\gamma)}+\frac{\mu^{\rho(1-\gamma)}}{\rho^{\alpha-\gamma}\Gamma(\alpha)}
\int_{0}^{t}\tau^{\rho-1}(t^{\rho}-\tau^{\rho})^{\alpha-1}\big[g(\tau) f((\tau^{\rho}/\rho)^{1-\gamma}|x(\tau)|)\big]d\tau
\\&\nonumber\quad\quad\quad\quad+\frac{\mu^{\rho(1-\gamma)}}{\rho^{\alpha-\gamma}\Gamma(\alpha)}
\int_{0}^{t}\tau^{\rho-1}(t^{\rho}-\tau^{\rho})^{\alpha-1}\psi(\tau)d\tau\\&\nonumber\leq
\frac{x_{0}}{\Gamma(\gamma)}+\frac{\mu^{\rho(1-\gamma)}}{\rho^{\alpha-\gamma}\Gamma(\alpha)}\|g\|_{[0,\mu]}
\int_{0}^{t}\tau^{\rho-1}(t^{\rho}-\tau^{\rho})^{\alpha-1}f((\tau^{\rho}/\rho)^{1-\gamma}|x(\tau)|)d\tau
\\&\quad\quad\quad\quad+\frac{\mu^{\rho(1-\gamma)}}{\rho^{\alpha-\gamma}\Gamma(\alpha)}
\int_{0}^{t}\tau^{\rho-1}(t^{\rho}-\tau^{\rho})^{\alpha-1}\psi(\tau)d\tau
\end{align}
Now, we taking 
\begin{align*}
\quad\phi(t)=\big(t^{\rho}/\rho\big)^{1-\gamma}|x(t)|,&\quad\theta(t)=\frac{x_{0}}{\Gamma(\gamma)}+
\frac{\mu^{\rho(1-\gamma)}}{\rho^{\alpha-\gamma}\Gamma(\alpha)}
\int_{0}^{t}\tau^{\rho-1}(t^{\rho}-\tau^{\rho})^{\alpha-1}\psi(\tau)d\tau~
and\\&\omega=\frac{\mu^{\rho(1-\gamma)}}{\rho^{\alpha-\gamma}\Gamma(\alpha)}\|g\|_{[0,\mu]}.
\end{align*} 
\par By applying Lemma 4.2, we can see that $\phi(t)=\big(t^{\rho}/\rho\big)^{1-\gamma}|x(t)|$ be a
bounded on $[0,\mu).$ Hence, for any $\zeta\in (0,\mu), x(t)$ be a
bounded on $[\zeta,\mu).$ By using Theorem 4.3, the initial value problems $(1.3)$ has a solution 
$x(t)~on~[0,+\infty).\quad\quad\quad\quad\quad\Box$
\par The next theorem guarantees the existence and uniqueness of global solution for the initial value problems $(1.3)~on~\mathbb{R}^{+}$
\\ \textbf{\ Theorem 4.4.} Suppose that $(\mathcal{H}_{1})$ is satisfied and there exists 
$g(t)$ be a non-negative continuous function defined on $[0,+\infty)$ such that
$\big|\varphi(t,x(t))-\varphi(t,y(t))\big|\leq g(t)|x(t)-y(t)|.$
Then the initial value problems $(1.3)$ has a unique solution in $C_{1-\gamma,\rho}[0,+\infty).$\\
\par We can obtained the existence of a global solution by using similar arguments as above. 
By applying Lipschitz condition and Lemma 4.2, we can deduced the uniqueness of global solution.
Here, we omitted the proof.
\section{Concluding Remarks}
\textbf{\ Remark 5.1.} If we take $\rho=1,$ in the initial value problems $(1.3)~and~(1.4),$ then
\par $(1)$ the local existence Theorems 3.1 and 3.2, yield the local existence Theorems 1 and 2, [21] respectively, 
associated with Hilfer-type fractional differential equations with the initial value problems $(5)~and~(6),$ [21] respectively.    
\par $(2)$ For $\beta=0,$ the local existence Theorems 3.1 yields the local existence Theorems 3.1 [20],
associated with Riemann–Liouville-type fractional differential equations with the initial value problems $(1)$ [21]. 
\par $(3)$ For $\beta=1$ the local existence Theorems 3.1 and 3.2, yield the local existence Theorems 3.1 and 3.2, [8] respectively,
associated with Caputo-type fractional differential equations with the initial value problems $(1.1)~and~(1.2),$ [8] respectively.
\\ \textbf{\ Remark 5.2.} If we take $\rho=1,$ in the initial value problems $(1.3)~and~(1.4),$ then
\par $(1)$ the continuation Theorems 4.1 and 4.2, reduce to the continuation Theorems 3 and 4, 
for Hilfer-type fractional differential equations [21] respectively.
\par $(2)$ For $\beta=0,$ the continuation Theorems 4.1 and 4.2, reduce to the continuation Theorems 4.1 and 4.2, 
for Riemann–Liouville-type fractional differential equations [20], respectively.
\par $(3)$ For $\beta=1$ the continuation Theorems 4.1 and 4.2, reduce to the continuation Theorems 4.2 and 4.4,
for Caputo-type fractional differential equations [8] respectively.
\\ \textbf{\ Remark 5.3.} In this article, we proved a new existence theorems of a local solutions for the generalized 
fractional differential equations, which is Hilfer-Katugampola-type with the certain singularity functions. 
Also, we obtained  two continuation theorems and we established global existence theorems for 
Hilfer-Katugampola-type fractional differential equations, which had been not investigated before. 
Our discussion in this article generalizes the existing results in the literature.

\end{document}